\documentclass[11pt]{amsart}
\usepackage{tabularx,booktabs}
\usepackage{caption}
\usepackage{amsmath}
\usepackage{amsfonts}

\usepackage{amscd}
\usepackage{amsthm}
\usepackage{amssymb} \usepackage{latexsym}
\usepackage{eufrak}
\usepackage{euscript}
\usepackage{epsfig}
\usepackage{graphics}
\usepackage{array}
\usepackage{enumerate}
\usepackage{dsfont}
\usepackage{color}
\usepackage{wasysym}
\usepackage{hyperref}
\usepackage{pdfsync}

\newcommand{\bel}[1]{\begin{equation}\label{#1}}

\newcommand{\be}{\begin{equation}}

\newcommand{\ba}{\begin{eqnarray}}
\newcommand{\ea}{\end{eqnarray}}

\newcommand{\qe}{\end{equation}}
\newcommand{\R}{{\mathbb R}}
\newcommand{\N}{{\mathbb N}}
\newcommand{\Z}{{\mathbb Z}}

\newcommand{\NNG}{\mathcal{PC}_{\geq 0}}

\newcommand{\Hmm}[1]{\leavevmode{\marginpar{\tiny%
$\hbox to 0mm{\hspace*{-0.5mm}$\leftarrow$\hss}%
\vcenter{\vrule depth 0.1mm height 0.1mm width \the\marginparwidth}%
\hbox to
0mm{\hss$\rightarrow$\hspace*{-0.5mm}}$\\\relax\raggedright #1}}}

\newtheorem{theorem}{Theorem}[section]

\newtheorem{lemma}[theorem]{Lemma}
\newtheorem{corollary}[theorem]{Corollary}

\newtheorem{remark}[theorem]{Remark}
\newtheorem{conjecture}[theorem]{Conjecture}
\newtheorem{prop}[theorem]{Proposition}

\begin{document}

\title[Total curvature of planar graphs with nonnegative curvature]{Total curvature of planar graphs with nonnegative combinatorial curvature}
\author{Bobo Hua}
\email{bobohua@fudan.edu.cn}
\address{School of Mathematical Sciences, LMNS, Fudan University, Shanghai 200433, China}

\author{Yanhui Su}
\email{suyh@fzu.edu.cn}
\address{College of Mathematics and Computer Science, Fuzhou University, Fuzhou 350116, China}

\begin{abstract} 
We prove that the total curvature of any planar graph with nonnegative combinatorial curvature is an integral multiple of $\frac{1}{12}.$ As a corollary, this answers a question proposed by T. R{\'e}ti.
\end{abstract}
\maketitle
\tableofcontents


\par
\maketitle

\bigskip


\section{Introduction}\label{s:intro}
The combinatorial curvature for a graph embedded into a surface was introduced by \cite{Neva70,Stone76,Gromov87,Ishida90}. The idea is to properly embed a graph into a piecewise flat surface and to use the generalized Gaussian curvature of the surface to define the combinatorial curvature for the graph. Many interesting geometric results have been obtained since then, see e.g. \cite{Zuk97,Woess98,Higuchi01,BP01,HJL02,LPZ02,HigS03,SunYu04,RBK05,BP06,DeVosMohar07,ChenChen08,Zhang08,Chenbf09,Keller10,KP11,Keller11,Oh17}.

Let $(V,E)$ be a (possibly infinite) locally finite, undirected simple graph with the set of vertices $V$ and the set of edges $E,$ which is topologically embedded into a surface $S$ without boundary. We write $G=(V,E,F)$ for the graph structure induced by the embedding where $F$ is the set of faces, i.e. connected components of the complement of the embedding image of the graph $(V,E)$ in $S,$ and call it a semiplanar graph, see \cite{HJL15} (It is called planar if $S$ is either the 2-sphere or the plane). We say that a semiplanar graph $G$ is a \emph{tessellation} of $S$ if the following hold, see e.g. \cite{Keller11}:
\begin{enumerate}[(i)]
\item Every face is homeomorphic to a disk whose boundary consists of finitely many edges of the graph.
\item Every edge is contained in exactly two different faces.
\item For any two faces whose closures have non-empty intersection, the intersection is either a vertex or an edge.
\end{enumerate} In this paper, we only consider tessellations of surfaces. 
For a semiplanar graph $G$, the \emph{combinatorial curvature} at the
vertex is defined as
\begin{equation}\label{def:comb}\Phi(x)=1-\frac{\deg(x)}{2}+\sum_{\sigma\in F:x\in \overline{\sigma}}\frac{1}{\deg(\sigma)},\quad x\in V,\end{equation} where the summation is taken over all faces $\sigma$ whose closure $\overline{\sigma}$ contains $x$ and $\deg(\cdot)$ denotes the degree of a vertex or a face, see Section~\ref{s:pre}. To digest the definition, we 
endow the ambient space $S$ with the metric structure: Each edge is of length one, each face is isometric to a Euclidean regular polygon of side length one with same facial degree which is glued along the edges, and the metric is induced by the gluing, see \cite[Definition~2.3.4]{BuragoBuragoIvanov01}. This piecewise flat surface is called the (regular) \emph{polygonal surface}, denoted by $S(G).$ It is well-known that the generalized Gaussian curvature, as a measure, concentrates on the vertices.
The combinatorial curvature at a vertex is in fact the mass of the generalized Gaussian curvature at that point up to the normalization $2\pi,$ see e.g. \cite{Alex05}. In this paper, we only consider the combinatorial curvature of a semiplanar graph and simply call it the curvature if it is clear in the context.

The Myers type theorem has been obtained by Stone \cite{Stone76}: A semiplanar graph with the curvature uniformly bounded below by a positive constant is a finite graph.
Higuchi \cite{Higuchi01} conjectured that it is finite even if the curvature is positive pointwise (i.e. at each vertex), which was proved by DeVos and Mohar 
\cite{DeVosMohar07}, see \cite{SunYu04} for the case of cubic graphs.  The largest size, i.e. the number of vertices, of planar graphs embedded into the 2-sphere, except the trivial classes including prisms and antiprisms, is known to be sandwiched between 208 and 380 by \cite{NichS11} and \cite{Oh17} respectively.

For a smooth surface with absolutely integrable Gaussian curvature, its total curvature encodes the global geometric information of the space, e.g. the boundary at infinity, see \cite{SST03}. For example, the total curvature of a convex surface in $\R^3$ describes the apex angle of the cone at infinity of the surface, which is useful to study global geometric and analytic properties of the surface, such as harmonic functions and heat kernels etc., following Colding and Minicozzi \cite{ColdingMinicozzi97JDG} and Xu \cite{Xu14}. In this paper, we study the total curvature of semiplanar graphs with nonnegative curvature. For a semiplanar graph $G,$ we denote by $$\Phi(G):=\sum_{x\in V}\Phi(x)$$ the total curvature of $G$ whenever the summation converges absolutely.  In case of finite graphs, the Gauss-Bonnet theorem reads as, see e.g. \cite{DeVosMohar07},
\begin{equation}\label{eq:GBfinite}\Phi(G)=\chi(S(G)),\end{equation} where $\chi(\cdot)$ denotes the Euler characteristic of a surface.
For an infinite semiplanar graph $G,$ the Cohn-Vossen type theorem, see \cite[Theorem~1.3]{DeVosMohar07} or \cite[Theorem~1.6]{ChenChen08}, yields that
\begin{equation}\label{CVthm}\Phi(G)\leq\chi(S(G)),\end{equation}
whenever $\sum_{x\in V}\min\{\Phi(x),0\}$ converges.

By the geometric meaning of the combinatorial curvature, one can see that a semiplanar graph $G$ has nonnegative combinatorial curvature if and only if the polygonal surface $S(G)$ is a generalized convex surface, i.e. an Alexandrov surface with nonnegative sectional curvature \cite{BuragoGromovPerelman92,BuragoBuragoIvanov01}. We denote by $\NNG$ the set of infinite semiplanar graphs with nonnegative combinatorial curvature. By the Cohn-Vossen type theorem, the inequality \eqref{CVthm} yields that $0\leq \Phi(G)\leq 1,$ for any $G\in \NNG.$ As is well-known in Riemannian geometry that for any real number $a\in[0,2\pi]$ there is a convex surface whose total curvature is given by $a.$ However, due to the combinatorial restriction of semiplanar graphs, by an elementary argument we prove that the set of all possible values of total curvatures of infinite semiplanar graphs with nonnegative curvature is a discrete subset in $[0,1],$ see Proposition \ref{prop:gap}. In particular, this yields that
$$\tau_1:=\inf\left\{\Phi(G):G\in \NNG,\Phi(G)>0\right\}>0,$$ which can be regarded as the first gap from zero for total curvatures of semiplanar graphs in the class of $\NNG.$ T. R\'eti, see \cite[Conjecture~2.1]{HuaLin16}, asked the question what $\tau_1$ is, and conjectured that $\tau_1=\frac{1}{6}$ and the minimum is attained by the graph consisting of a pentagon and infinitely many hexagons. 

\begin{figure}[htbp]
\begin{center}
\includegraphics[width=0.8\linewidth]{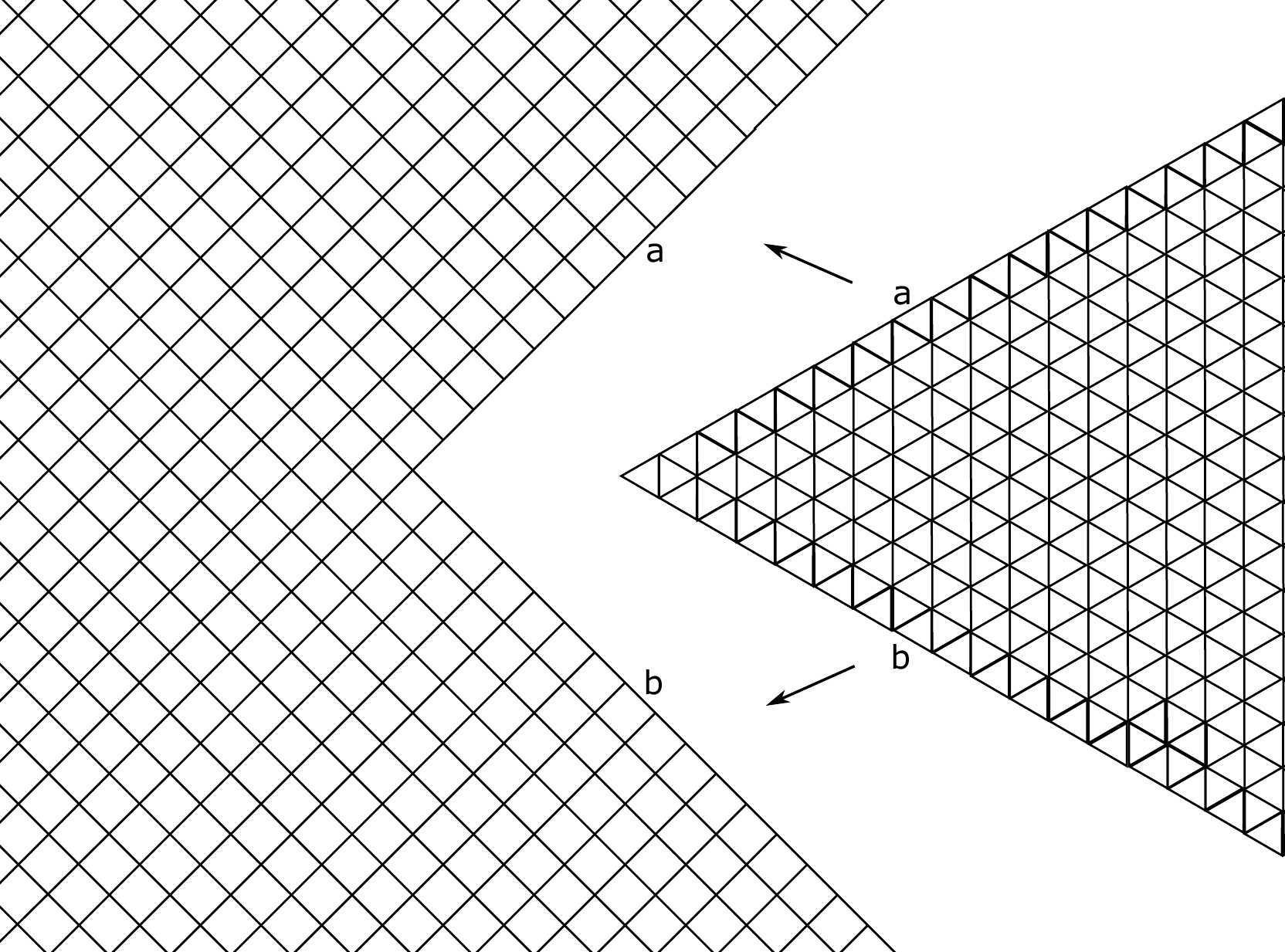}
\caption{\small A graph with total curvature $\frac{1}{12}$ in which half lines with same labeling are glued together.}
\label{fig2}
\end{center}
\end{figure}

\begin{figure}[tb]
\begin{minipage}[b]{\textwidth}
\centering
\includegraphics[width=6cm,height=7.7cm]{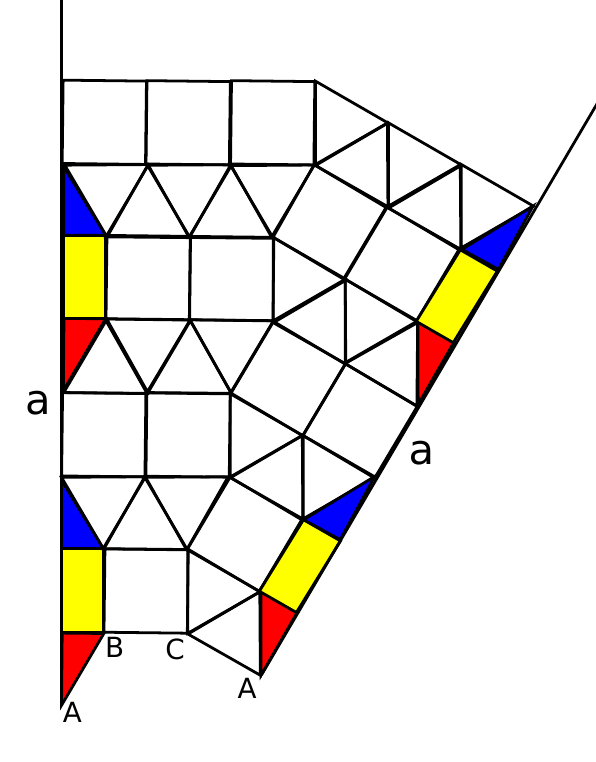}
\end{minipage}
\caption{\small A graph with total curvature $\frac{11}{12}$ obtained from gluing the sides $a$ in the picture and attaching another triangle to the vertices $A,B$ and $C.$
}
 \label{11d12-1}
\end{figure}

Inspired by R\'eti's question, it will be interesting to know all possible values of total curvatures of semiplanar graphs with nonnegative curvature. Using the combinatorial information and the Gauss-Bonnet theorem (for compact subsets  with boundary), we are able to determine all values of total curvatures in this class.
\begin{theorem}\label{mainthm2}
The set of all values of total curvatures of infinite semiplanar graphs with nonnegative combinatorial curvature is given by $$\left\{\frac{i}{12}:\ 0\leq i\leq 12, i\in \Z\right\}.$$
\end{theorem} As a corollary, this answers R\'eti's question that $\tau_1=\frac{1}{12}$ and the minimum can be attained by e.g. Figure~\ref{fig2}. For any figure in this paper, vertices, edges or sides on the boundary with same labeling are identified with each other.  While it is hard to classify the graph/tessellation structures for semiplanar graphs which attain the first gap of the total curvature, see \cite{GrunbaumS87} for classification problems of planar tilings with vanishing curvature, the metric structures for their ambient polygonal surfaces can be classified, see a companion paper \cite{HuaSu17}.
Moreover, as the part of the theorem, we construct semiplanar graphs with nonnegative curvature whose total curvatures attain all values listed above, see Figure~\ref{11d12-1} for a semiplanar graph with total curvature $\frac{11}{12}$ and other examples in Section~\ref{s:examples}. 

These proof strategies apply to the problems on the total curvature of a semiplanar graph with boundary, i.e. a graph embedded into a surface with boundary. For precise definitions and the terminology used, we refer to Section~\ref{s:boundary}. Let $G$ be a semiplanar graph with boundary and with nonnegative curvature, and $S(G)$ be the polygonal surface of $G.$ Consider the doubling constructions of $S(G)$ and $G,$ see e.g. \cite[Section~5]{Perelman91}. Let $\widetilde{S(G)}$ be the
double of ${S}(G)$, that is, $\widetilde{S(G)}$ consists of two copies of ${S}(G)$ glued along the boundaries via the identity map restricted on $\partial S(G)$.
This induces the doubling graph of $G,$ denoted by $\widetilde{G}.$ By the definition of the curvature for semiplanar graphs with boundary, one can show that $\widetilde{G}$ has nonnegative curvature and $\Phi(\widetilde{G})=2\Phi(G).$ This yields that the total curvature of $G$ is an integral multiple of $\frac{1}{24}$ by Theorem~\ref{mainthm2}. In fact, we can prove stronger results for semiplanar graphs with boundary:
\begin{itemize}
\item The total curvature of a semiplanar graph with boundary and with nonnegative curvature is an integral multiple of $\frac{1}{12},$ see Theorem~\ref{half-plane-2}.
\item The graph/tessellation structures for the semigraphs with boundary attaining the first gap of the total curvature, $\frac{1}{12}$ by the above result, are classified in Theorem~\ref{thm:classification}.
\end{itemize}

The paper is organized as follows: In next section, we recall some basics on the combinatorial curvature for semiplanar graphs.
Section~\ref{s:thm2} is devoted to the proof of Theorem~\ref{mainthm2}. 
In Section~\ref{s:examples}, we construct examples to show the sharpness of Theorem~\ref{mainthm2}. In the last section, we study the total curvature for semiplanar graphs with boundary.


\section{Preliminaries}\label{s:pre}
Let $G=(V,E,F)$ be a semiplanar graph induced by an embedding of a graph $(V,E)$ into a (possibly un-orientable) surface $S$ without boundary.  We only consider the appropriate embedding such that $G$ is a tessellation of $S,$ see the definition in the introduction.

We say that a vertex is incident to an edge (similarly, an edge is incident to a face, or a vertex is incident to a face) if the former is a subset of the closure of the latter. Two vertices are called neighbors if there is an edge connecting them. We denote by $\deg(x)$ the number of neighbors of a vertex $x,$ and by $\deg(\sigma)$ the number of edges incident to a face $\sigma$ (equivalently, the number of vertices incident to $\sigma$). Two edges (two faces resp.) are called adjacent if there is a vertex (an edge resp.) incident to both of them.
The combinatorial distance between two vertices $x$ and $y,$ denote by $d(x,y),$ is defined as the minimal length of walks from $x$ to $y,$ i.e. the minimal number $n$ such that there is $\{x_{i}\}_{i=1}^{n-1}\subset V$ satisfying $x\sim x_1\sim\cdots\sim x_{n-1}\sim y.$ We denote by $B_r(x):=\{y\in V: d(y,x)\leq r\},$ $r\geq 0,$ the ball of radius $r$ centered at the vertex $x.$ 
For a tessellation, we always assume that $3\leq \deg(x)<\infty$ and $3\leq\mathrm{deg}(\sigma)<\infty$ for any vertex $x$ and face $\sigma.$

Given a semiplanar graph $G=(V,E,F)$ embedded into a surface $S,$ it associates with a unique metric space $S(G),$ called the polygonal surface, which is obtained from replacing each face of $G$ by a regular Euclidean polygon of side length one with same facial degree and gluing the polygons along the edges. Note that the (induced) metric on $S(G)$ is piecewise flat and generally non-smooth near a vertex, while it is locally isometric to a flat domain in $\R^2$ near any interior point of an edge or a face. As a metric surface with isolated singularities, the generalized Gaussian curvature $K$ of $S(G)$ vanishes at smooth points and can be regarded as a measure concentrated on the singularities, i.e. on vertices. One can show that the mass of the generalized Gaussian curvature at each vertex $x$ is given by
$K(x) = 2\pi-\Sigma_x,$
where $\Sigma_x$ denotes the total angle at $x$ in the metric space $S(G),$ see \cite{Alex05}. Moreover, by direct computation one has
$K(x) = 2\pi\Phi(x),$ where the combinatorial curvature $\Phi(x)$ is defined in \eqref{def:comb}, namely, the combinatorial curvature at each vertex is equal to the mass of the generalized Gaussian curvature at that point up to the normalization $2\pi.$ Hence one can show that a semiplanar graph $G$ has nonnegative combinatorial curvature if and only if the polygonal surface $S(G)$ is a convex surface.

In this paper, we consider the total curvature of infinite semiplanar graphs with nonnegative combinatorial curvature. For our purposes, it suffices to consider those with positive total curvature (otherwise the total curvature vanishes). By \cite[Theorem 3.10]{HJL15}, these graphs are planar, namely, the ambient spaces are homeomorphic to $\R^2.$

In a semiplanar graph, a pattern of a vertex $x$ is defined as a vector $(\mathrm{deg}(\sigma_1),\mathrm{deg}(\sigma_2),\cdots,\mathrm{deg}(\sigma_{N})),$ where $N=\deg(x),$ $\{\sigma_i\}_{i=1}^N$ are the faces which $x$ is incident to, and $\mathrm{deg}(\sigma_1)\leq\mathrm{deg}(\sigma_2)\leq\cdots\leq\mathrm{deg}(\sigma_{N}).$ For simplicity, we always write $$x=(\mathrm{deg}(\sigma_1),\mathrm{deg}(\sigma_2),\cdots,\mathrm{deg}(\sigma_{N}))$$ to indicate the pattern of the vertex $x.$

Table \ref{tabl1} is the list of all
possible patterns of a vertex with positive curvature (see
\cite{DeVosMohar07,ChenChen08}); Table \ref{tabl2} is the list of all
possible patterns of a vertex with vanishing curvature (see \cite{GrunbaumS87,ChenChen08}).

\begin{table}
\refstepcounter{table}\label{tabl1}
\begin{tabular}{|lc|lr|}
\hline
Patterns &&$\Phi(x)$&\\
    \hline
 $(3,3,k)$ & $3\leq k$&$1/6+1/k$&\\
 $(3,4,k)$  & $4\leq k$  &$1/12+1/k$&\\
 $(3,5,k)$ & $5\leq k$&$1/30+1/k$&\\
$(3,6,k)$&$6\leq k$&$1/k$&\\
$(3,7,k)$&$7\leq k\leq41$&$1/k-1/42$&\\
$(3,8,k)$&$8\leq k\leq 23$&$1/k-1/24$&\\
$(3,9,k)$&$9\leq k\leq 17$&$1/k-1/18$&\\
$(3,10,k)$&$10\leq k\leq 14$&$1/k-1/15$&\\
$(3,11,k)$&$11\leq k\leq 13$&$1/k-5/66$&\\
$(4,4,k)$&$4\leq k$&$1/k$&\\
$(4,5,k)$&$5\leq k\leq 19$&$1/k-1/20$&\\
$(4,6,k)$&$6\leq k\leq 11$&$1/k-1/12$&\\
$(4,7,k)$&$7\leq k\leq 9$&$1/k-3/28$&\\
$(5,5,k)$&$5\leq k\leq 9$&$1/k-1/10$&\\
$(5,6,k)$&$6\leq k\leq 7$&$1/k-2/15$&\\
$(3,3,3,k)$&$3\leq k$&$1/k$&\\
$(3,3,4,k)$&$4\leq k\leq 11$&$1/k-1/12$&\\
$(3,3,5,k)$&$5\leq k\leq 7$&$1/k-2/15$&\\
$(3,4,4,k)$&$4\leq k\leq 5$&$1/k-1/6$&\\
$(3,3,3,3,k)$&$3\leq k\leq5$&$1/k-1/6$&\\
\hline
\multicolumn{4}{c}{}\\    
\end{tabular}

\textbf{\tablename~\thetable.} The patterns of a vertex with positive curvature.
\end{table}

\begin{table}
\refstepcounter{table}\label{tabl2}
\begin{tabular}{lllll}
\hline
$(3,7,42),$ & $(3,8,24),$ & $(3,9,18),$ & $(3,10,15),$ & $(3,12,12),$\\
$(4,5,20),$ & $(4,6,12),$ & $(4,8,8),$ & $(5,5,10),$ & $(6,6,6),$\\
$(3,3,4,12),$ & $(3,3,6,6),$& $(3,4,4,6),$ & $(4,4,4,4),$ &$(3,3,3,3,6),$\\
$(3,3,3,4,4),$&$(3,3,3,3,3,3).$&&&\\
\hline
\multicolumn{1}{c}{}\\    
\end{tabular}

\textbf{\tablename~\thetable.} The patterns of a vertex with vanishing curvature.
\end{table}

%
%

For any infinite semiplanar graph with nonnegative curvature, Chen and Chen \cite{ChenChen08,Chenbf09} obtained
an interesting result that the curvature vanishes outside a finite subset of vertices in the graph.
\begin{theorem}[Theorem~1.4 in \cite{ChenChen08}, Theorem~3.5 in \cite{Chenbf09}]\label{thm:ccthm} For a semiplanar graph with nonnegative curvature, there are only finitely
many vertices with non-vanishing curvature.
\end{theorem}
The following lemma is useful in this paper, see \cite[Lemma~2.5]{ChenChen08}.
\begin{lemma}\label{lemma}
If there is a face $\sigma$ such that $\mathrm{deg}(\sigma)\geq43$ and $\Phi(x)\geq0$ for
any vertex $x$ incident to $\sigma,$ then
$$\sum_{x}\Phi(x)\geq1,$$ where the summation is taken over all vertices $x$ incident to $\sigma.$
\end{lemma}

In fact, for an infinite semiplanar graph $G$ with nonnegative curvature if there is a face whose degree is at least $43,$ then the graph has rather special
structure, see \cite[Theorem~2.10]{HJL15}. In particular, the total curvature of $G$ is equal to $1$ and the polygonal surface $S(G)$ is isometric to a half flat-cylinder
in $\mathbb{R}^3.$
%
%
%




For our purposes, we need the Gauss-Bonnet theorem for piecewise flat manifolds with boundary, see e.g. \cite[Theorem~3.1]{DeVosMohar07}. Let $K\subset S(G)$ be the closure of a finite union of faces, regarded as a manifold with boundary. Then 
\begin{equation}\label{eq:GaussBonnetboundary}2\pi\sum_{x\in V\cap \mathrm{int}({K})}\Phi(x)+\sum_{x\in V\cap \partial K}(\pi-\theta_x)=2\pi\chi(K),\end{equation} where $\mathrm{int(\cdot)}$ denotes the interior of a set and $\theta_x$ is the inner angle at the boundary vertex $x$ w.r.t. the set $K.$

Since all possible patterns of a vertex with positive curvature are listed in Table \ref{tabl1}, we adopt a number-theoretic argument to show the restriction of the total curvature in the class of $\NNG,$ i.e. infinite semiplanar graphs with nonnegative curvature.
\begin{prop}\label{prop:gap}
The total curvature of a semiplanar graph with nonnegative curvature is an integral multiple of $1/219060189739591200.$ In particular, the set of all values of total curvatures for graphs in the class of $\NNG$ is discrete in $[0,1].$
\end{prop} 

\begin{proof} By the Gauss-Bonnet theorem for finite graphs \eqref{eq:GBfinite}, it suffices to consider infinite graphs. By Lemma \ref{lemma} and \eqref{CVthm}, or \cite[Theorem~2.10]{HJL15}, we only need to consider $G\in\mathcal{PC}_{\geq0}$ with the maximum facial degree less than or equal to $42,$ otherwise the total curvature is equal to $1.$ In this case, all patterns of a vertex with positive curvature are shown in Table \ref{tabl1} and the total curvature of $G$ is the sum of the curvatures at finitely many vertices with non-vanishing curvature. Thus, the total curvature is an integral multiple of $1/219060189739591200,$ where the denominator is given by
$$2^5\times3^3\times5^2\times7\times11\times13\times17\times19\times23\times29\times31\times37\times41$$ which
is the least common multiple of the denominators of reduced fractions of all possible values of positive curvature shown in Table \ref{tabl1}.
\end{proof}

\section{Proof of Theorem \ref{mainthm2}}\label{s:thm2}
Let $G$ be an infinite planar graph with nonnegative curvature and $S(G)$ be the polygonal surface of $G$. We denote by $T(G):=\{x\in V: \Phi(x)\neq 0\}$ the set of vertices of $G$ with non-vanishing curvature, which is a finite set by \cite{ChenChen08,Chenbf09}, see Theorem~\ref{thm:ccthm} in Section~\ref{s:pre}. 
Given a vertex $x$ with vanishing curvature, its possible patterns are listed in Table \ref{tabl2}. We may further reduce the pattern list if the curvature vanishes on $B_2(x).$

\begin{figure}[tb]
\begin{minipage}[b]{0.49\textwidth}
\centering
\includegraphics[width=5.1cm,height=3.47cm]{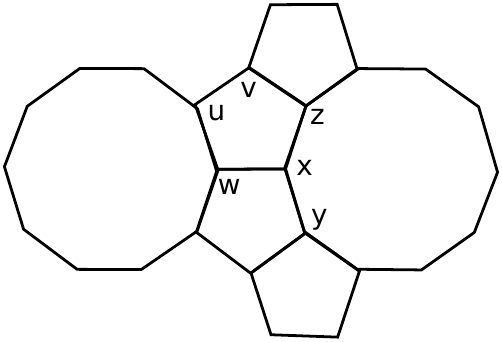}
\end{minipage}
\caption{\small Proof of Lemma \ref{lemma2.1}
}
 \label{5-5-10}
\end{figure}
\begin{lemma}\label{lemma2.1}
For any vertex $x$ satisfying $\Phi(y)=0$ for any $y\in B_2(x),$ the possible patterns of $x$ are in the following:
$$(3,12,12), (4,6,12), (4,8,8), (6,6,6), (3,3,4,12), (3,3,6,6), (3,4,4,6),$$
$$(4,4,4,4), (3,3,3,3,6), (3,3,3,4,4), (3,3,3,3,3,3).$$
\end{lemma}
\begin{proof} It suffices to exclude the following patterns in Table \ref{tabl2}, $$(3,7,42),(3,8,24),(3,9,18),(3,10,15),(4,5,20),(5,5,10).$$ 
For example, we consider the case of $x=(5,5,10).$ For other cases, the proofs are similar and hence omitted here. Let $y,z$ and $w$ be the neighbors of $x,$ see Figure \ref{5-5-10}. By the assumption, $y,z$ and $w$ have vanishing curvature, which yields that they are of the pattern $(5,5,10).$ For vertices $u$ and $v$ as in Figure \ref{5-5-10}, $\Phi(u)=\Phi(v)=0$ by the assumption. For the vertex $u,$ we have $u=(5,5,10),$ so that $v$ is incident to at least three pentagons and hence has non-vanishing curvature. This yields a contradiction and we can rule out the case of $x=(5,5,10).$
\end{proof}

For our purposes, let $G$ be an infinite planar graph with nonnegative curvature. Since $T(G)$ is a finite set and the polygonal surface $S(G)$ is homeomorphic to $\R^2,$ one can find a (sufficient large) compact set $K\subset S(G)$ satisfying the following:
\begin{itemize}
\item $K$ is the closure of a finite union of faces in $S(G)$ and is simply connected.
\item $T(G)\subset K.$
\item The boundary of $K,$ denoted by $\partial K,$ is a Jordan curve consisting of finitely many edges.
\item For any vertex $x\in\partial K,$ $B_2(x)\cap T(G)=\emptyset.$
\end{itemize}

\begin{lemma}\label{lem:488} Let $K$ satisfy the above properties. Suppose that there exists a vertex on $\partial K$ with the pattern $(4,8,8),$ then all patterns of vertices on $\partial K$ are $(4,8,8).$ 
\end{lemma}
\begin{proof} For any $x\in \partial K\cap V$ with the pattern $(4,8,8),$ Lemma \ref{lemma2.1} yields that for any neighbor vertex of $x$ on $\partial K$ has the same pattern $(4,8,8).$ Then the lemma follows from the connectedness of $\partial K.$
\end{proof}
\begin{remark} In above case, by the same argument, one can show that any vertex outside $K$ also has the pattern $(4,8,8).$
\end{remark}


\begin{proof}[Proof of Theorem \ref{mainthm2}]
For the first part of the theorem, we prove that all possible values of total curvatures of infinite planar graphs with nonnegative curvature are $\frac{i}{12}$, ($i\in\mathbb{Z}$ and $0\leq i\leq12$).  We apply the Gauss-Bonnet theorem \eqref{eq:GaussBonnetboundary} to $K,$ 
$$2\pi\sum_{x\in V\cap \mathrm{int}(K)}\Phi(x)+\sum_k(\pi-\theta_k)=2\pi,$$
where $\theta_k$ are the inner angles of vertices on $\partial K.$ By the properties of $K,$ we have $\Phi(G)=\sum_{x\in V\cap \mathrm{int}(K)}\Phi(x).$
Hence it suffices to prove that $\sum_k(\pi-\theta_k)$ is an integral multiple of $\frac{\pi}{6}$.
We divided into two cases.
\begin{description}
\item[Case 1] There is a vertex on $\partial K$ of the pattern $(4,8,8).$ By Lemma~\ref{lem:488}, all the vertices on $\partial K$ have the same pattern. Since $\partial K$ consists of edges, we divide them into two classes: \\
$\alpha$:=\{Edges\ in\ $\partial K$\ which\ are\ incident to a square and an octagon\} and
 $\beta$:=\{Edges\ in\ $\partial K$\ which\ are\ incident to two octagons\}.
Note that two distinct edges in the class $\beta$ are not adjacent by the combinatorial restriction of tilings by $(4,8,8).$ Since $\partial K$ is connected, the class $\alpha$ is non-empty. Take a closed walk along $\partial K$ starting from an edge $e_1\in\alpha$ and ending at the same edge without repeated edges, say,
\begin{equation}\label{eq:equwalk}e_1 f_1e_2 e_3f_2 e_4f_3e_5\cdots e_N f_Le_1,\end{equation} where $\{e_i\}_{i=1}^N\subset \alpha, \{f_j\}_{j=1}^L\subset \beta$ for some $N,L\in \N,$ and two successive edges are adjacent. For convenience, we write $e_{N+1}=e_1.$ For each $f_j\in \beta,$ $1\leq j\leq L,$ its adjacent edges in the sequence \eqref{eq:equwalk} are in the class $\alpha.$ For any two adjacent edges in \eqref{eq:equwalk}, say $e$ and $f,$ on $\partial K,$ we denote by $\theta_{ef}$ its inner angle at the common end-vertex of $e$ and $f$ w.r.t. the set $K.$ By the tiling $(4,8,8),$ for any $e_i,e_{i+1}\in \alpha$ and $f_j\in \beta,$ $$\theta_{e_i e_{i+1}}=\frac{\pi}{2}\ \mathrm{or}\ \frac{3\pi}{2},\quad \theta_{e_i f_j}=\frac{3\pi}{4}\ \mathrm{or}\ \frac{5\pi}{4}.$$ For any $1\leq i\leq N,$ we define
$\widehat{\theta_{i}}:=\theta_{e_i f_j}+\theta_{f_j e_{i+1}}$ if there is an edge $f_j\in \beta$ which is adjacent to  $e_i$ and $e_{i+1}$ in the walk \eqref{eq:equwalk} and $\widehat{\theta_{i}}:= \theta_{e_ie_{i+1}},$ otherwise. This yields that $\widehat{\theta_i}=\frac{\pi}{2},\frac{3\pi}{2}, \frac{5\pi}{2}\ \mathrm{or}\ 2\pi, \forall\ 1\leq i\leq N,$ which are integral multiples of $\frac{\pi}{2}.$ Hence
$$\sum_k(\pi-\theta_k)=M \pi-\sum_{i=1}^N \widehat{\theta_i},$$ where $M$ is the number of vertices on $\partial K,$ which implies that $\sum_k(\pi-\theta_k)$ is an integral multiple of $\frac{\pi}{2}.$

\item[Case 2] There is no vertices of the pattern $(4,8,8)$ on $\partial K.$ 
\begin{figure}[tb]
\begin{minipage}[b]{0.49\textwidth}
\centering
\includegraphics[width=2.967cm,height=2.967cm]{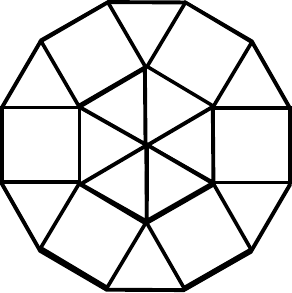}
\end{minipage}
\caption{\small A regular dodecagon can be divided into triangles and squares.
}
\label{12=34}
\end{figure}
Note that a regular dodecagon can be divided into twelve triangles and six squares, see Figure \ref{12=34}; a regular hexagon can be divided into six triangles. By dividing dodecagons and hexagons into triangles and squares, without loss of generality, we may assume that all vertices on $\partial K$ are incident to triangles or squares by Lemma~\ref{lemma2.1}. Thus the inner angle of any vertex on $\partial K$ are of the form $m\frac{\pi}{2}+n\frac{\pi}{3}$ $(m,n\in\mathbb{Z})$, which yields that $\sum_k(\pi-\theta_k)$ is an integral multiple of $\frac{\pi}{6}$.
\end{description}

Combining these results, we prove the first part of theorem \ref{mainthm2}.

For the other part of the theorem, we construct several graphs with nonnegative curvature whose total curvatures attain the full set $$\left\{\frac{i}{12}:\ 0\leq i\leq 12, i\in \Z\right\},$$ see Section~\ref{s:examples}.
\end{proof}

\section{Examples}\label{s:examples}

\begin{figure}[htbp]
\begin{center}
\includegraphics[width=0.4\linewidth]{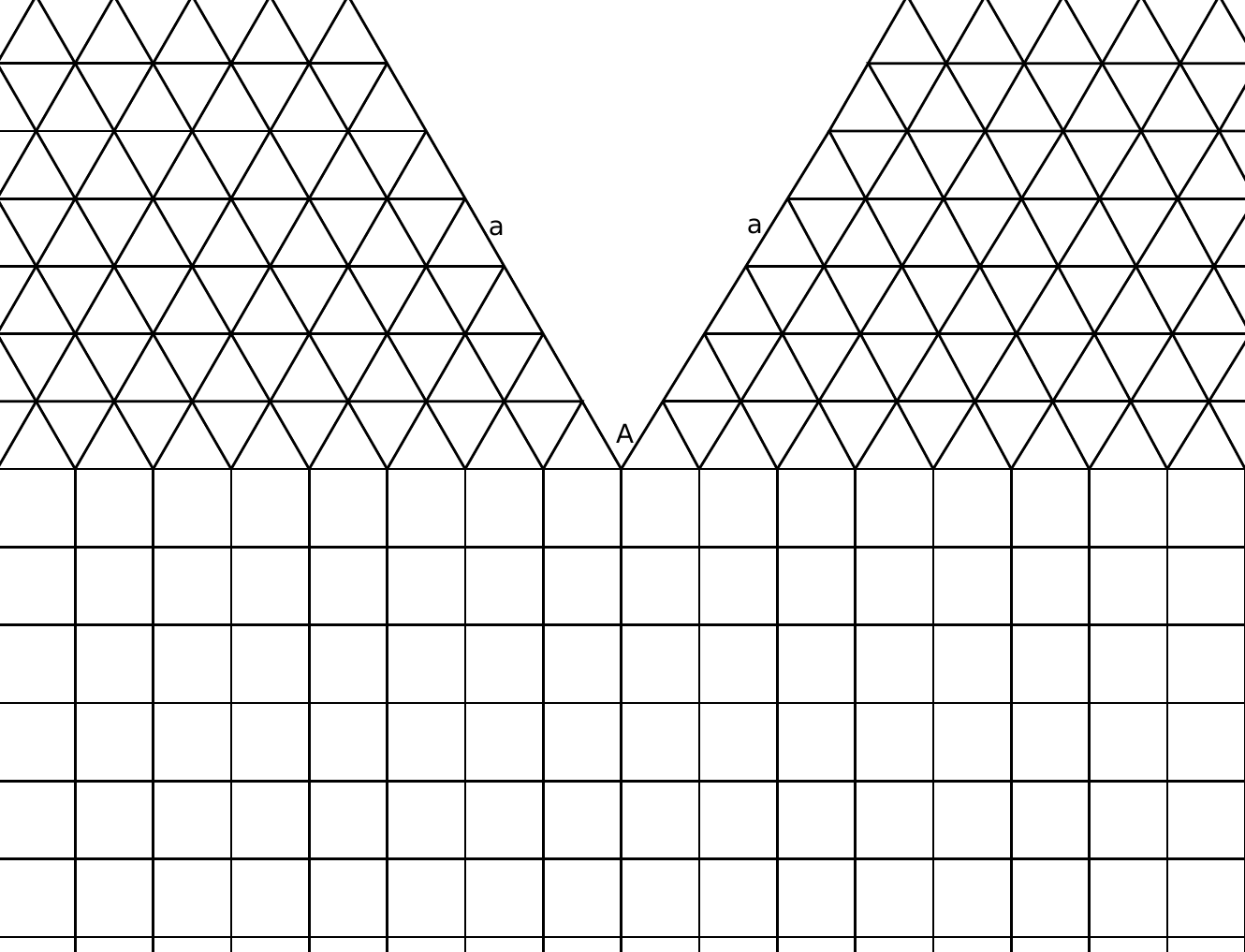}
\caption{\small A graph in $\mathcal{PC}_{\geq0}$ with total curvature $\frac{2}{12}.$}
\label{2d12}
\end{center}
\end{figure}

\begin{figure}[htbp]
\begin{center}
\includegraphics[width=0.4\linewidth]{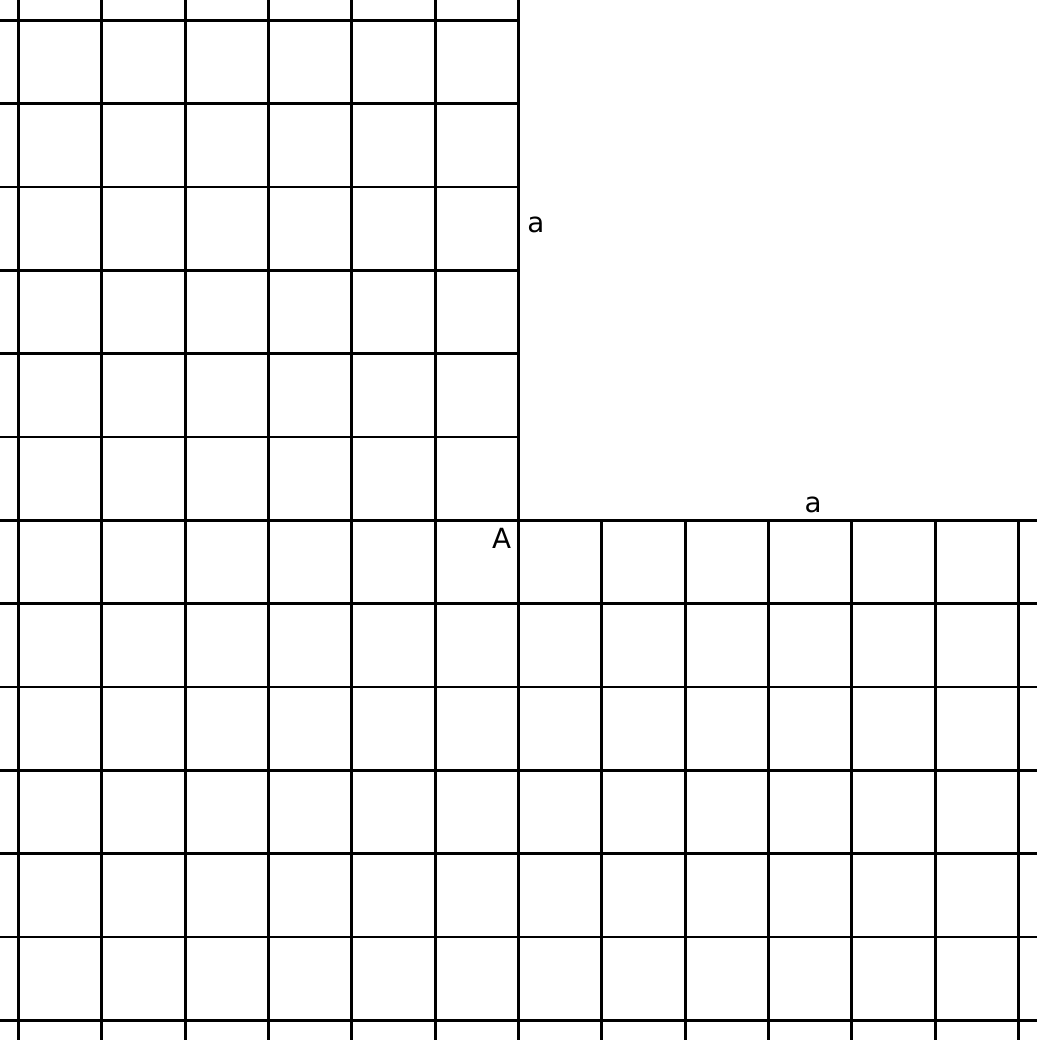}
\caption{\small A graph in $\mathcal{PC}_{\geq0}$ with total curvature $\frac{3}{12}$.}
\label{3d12}
\end{center}
\end{figure}

\begin{figure}[htbp]
\begin{center}
\includegraphics[width=0.4\linewidth]{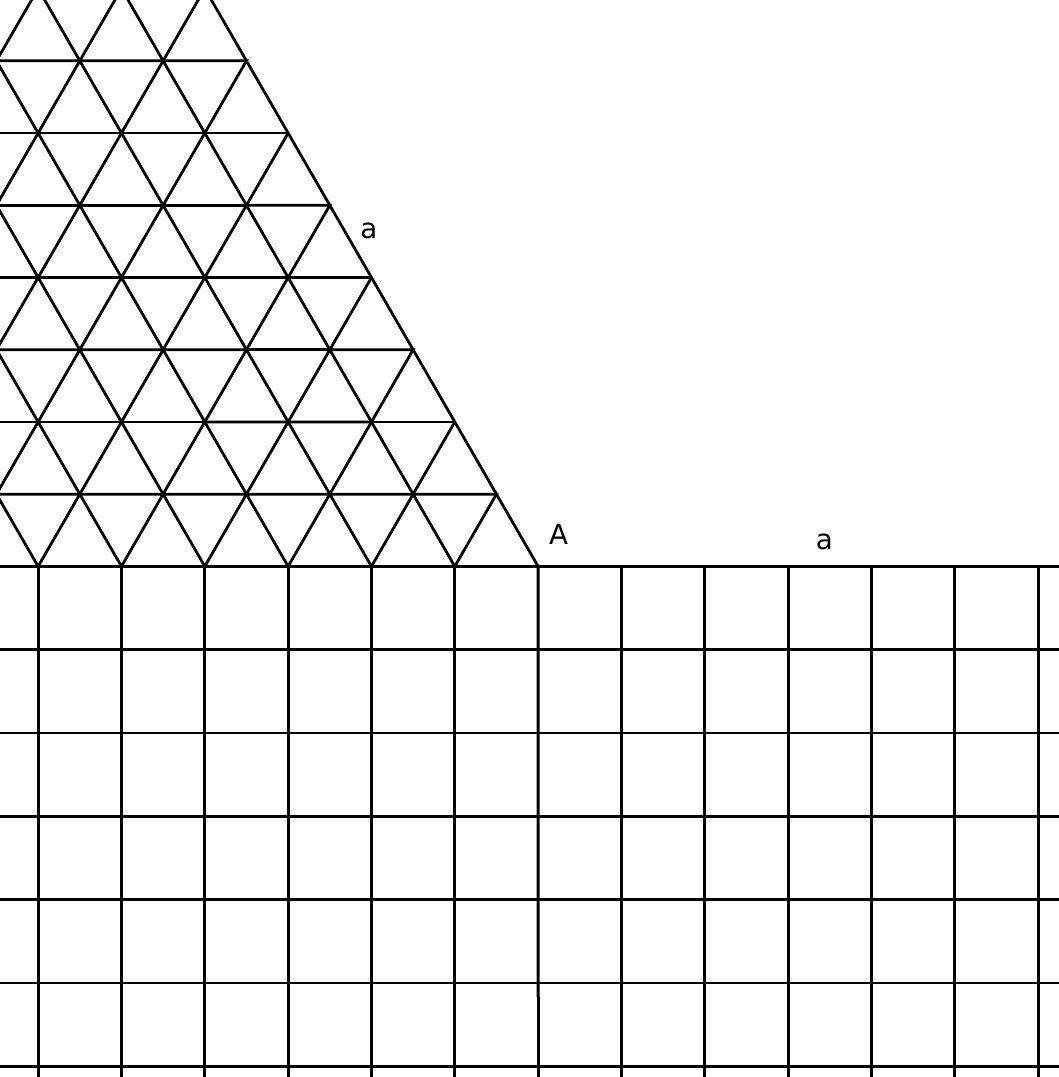}
\caption{\small A graph in $\mathcal{PC}_{\geq0}$ with total curvature $\frac{4}{12}$.}
\label{4d12}
\end{center}
\end{figure}

\begin{figure}[htbp]
\begin{center}
\includegraphics[width=0.4\linewidth]{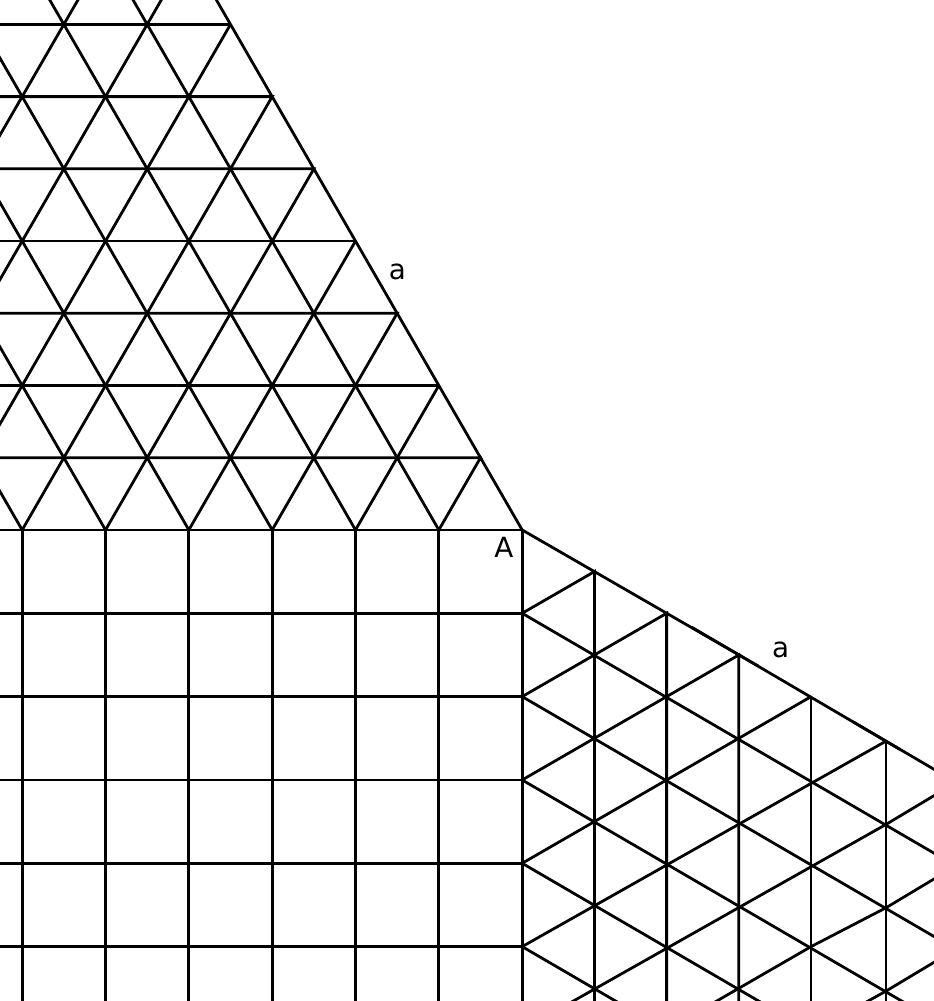}
\caption{\small A graph in $\mathcal{PC}_{\geq0}$ with total curvature $\frac{5}{12}$.}
\label{5d12}
\end{center}
\end{figure}

\begin{figure}[htbp]
\begin{center}
\includegraphics[width=0.4\linewidth]{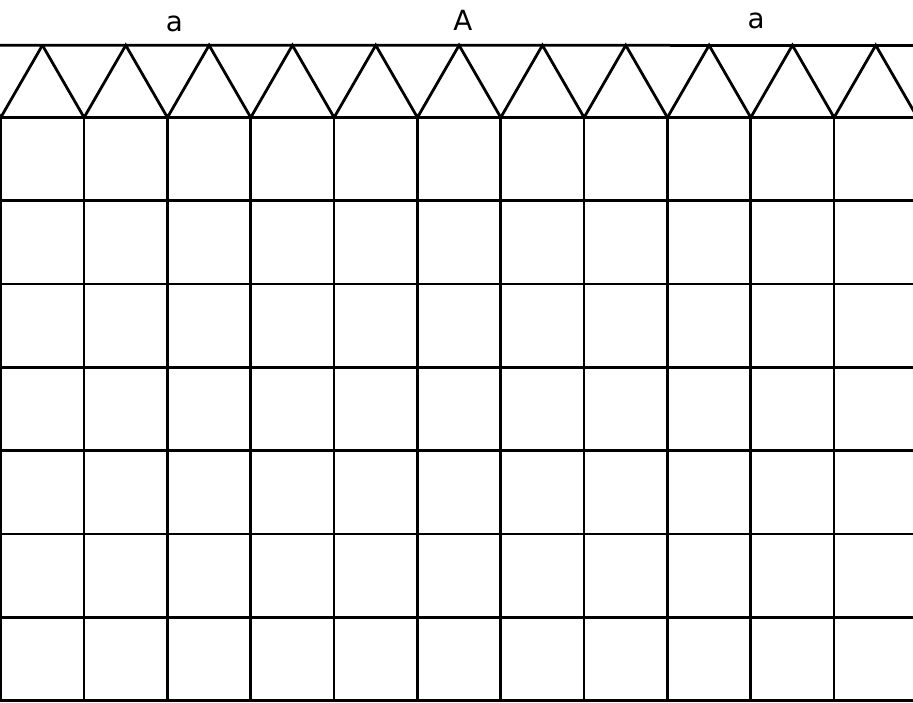}
\caption{\small A graph in $\mathcal{PC}_{\geq0}$ with total curvature $\frac{6}{12}$.}
\label{6d12}
\end{center}
\end{figure}

\begin{figure}[htbp]
\begin{center}
\includegraphics[width=0.4\linewidth]{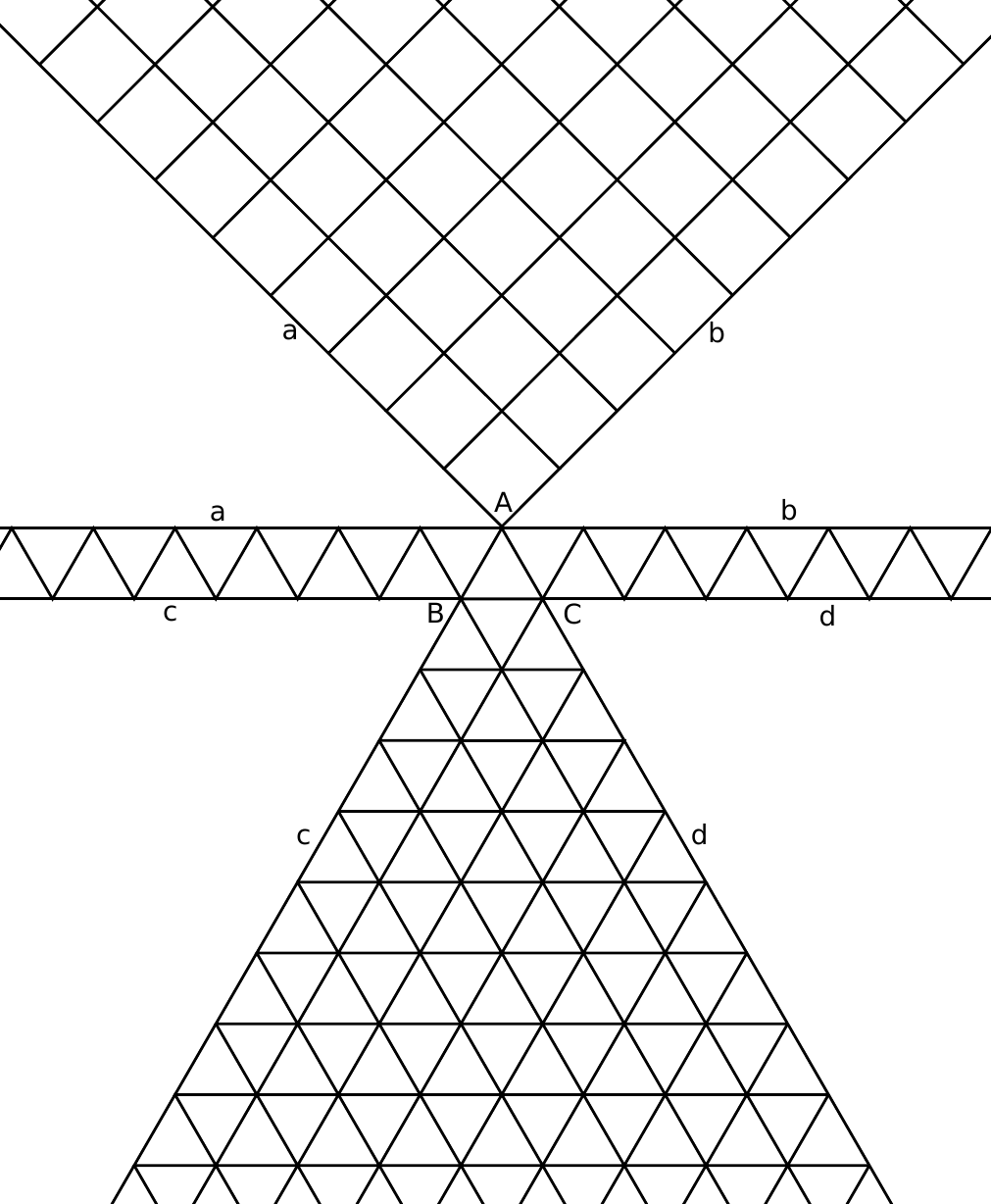}
\caption{\small A graph in $\mathcal{PC}_{\geq0}$ with total curvature $\frac{7}{12}.$}
\label{7d12}
\end{center}
\end{figure}

\begin{figure}[htbp]
\begin{center}
\includegraphics[width=0.4\linewidth]{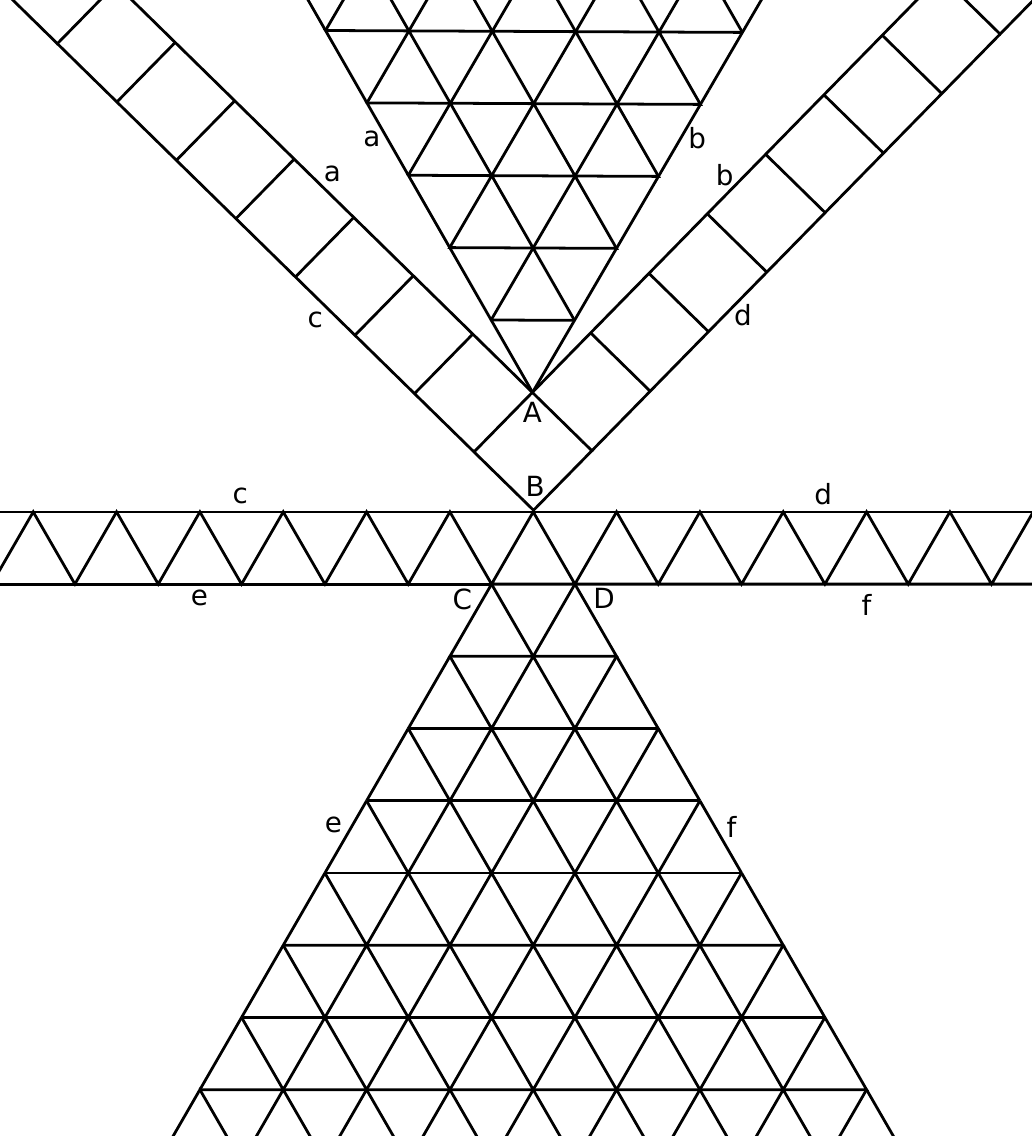}
\caption{\small A graph in $\mathcal{PC}_{\geq0}$ with total curvature $\frac{8}{12}.$}
\label{8d12}
\end{center}
\end{figure}

\begin{figure}[htbp]
\begin{center}
\includegraphics[width=0.5\linewidth]{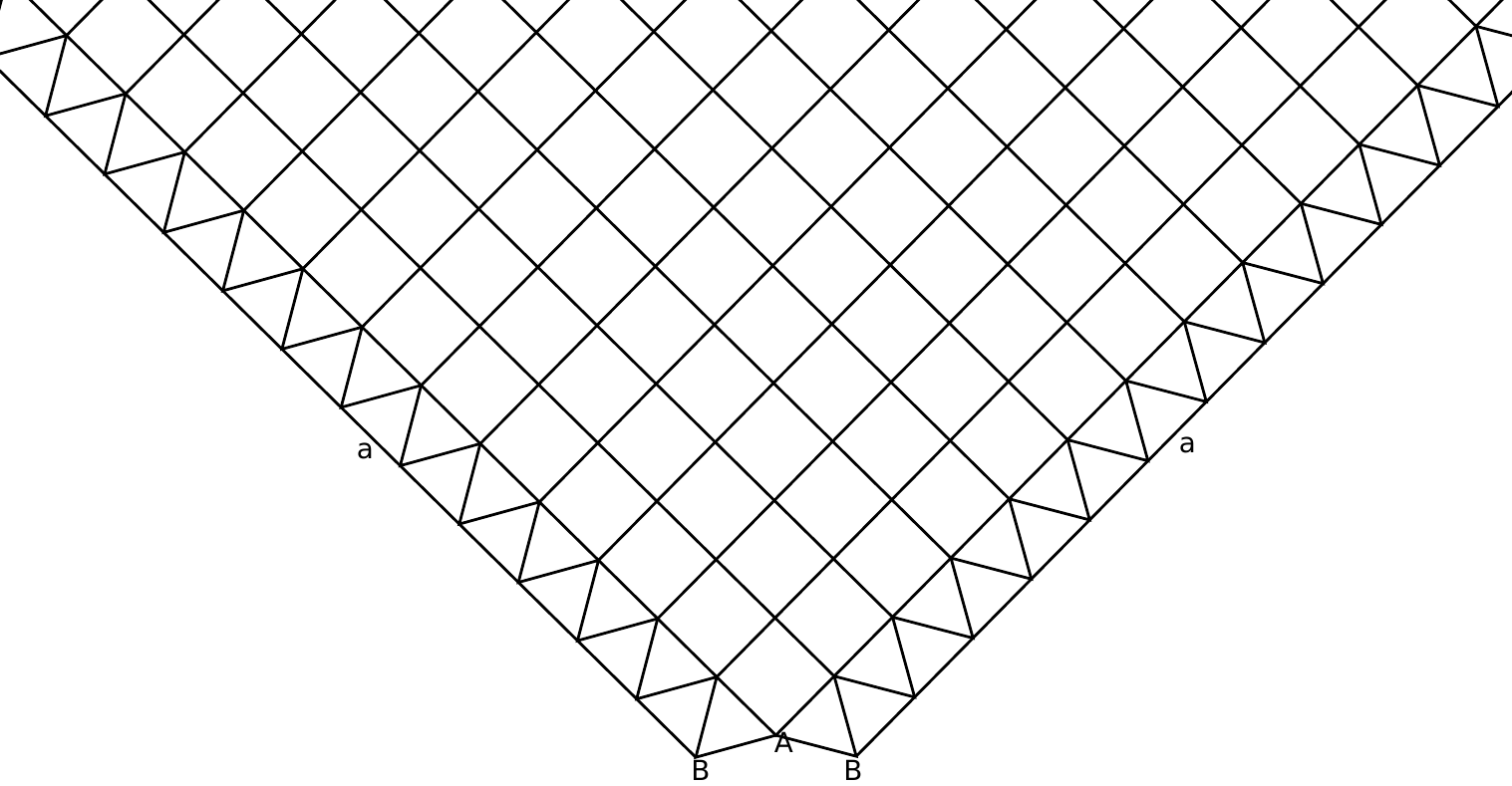}
\caption{\small A graph in $\mathcal{PC}_{\geq0}$ with total curvature $\frac{9}{12}$.}
\label{9d12}
\end{center}
\end{figure}

\begin{figure}[htbp]
\begin{center}
\includegraphics[width=0.4\linewidth]{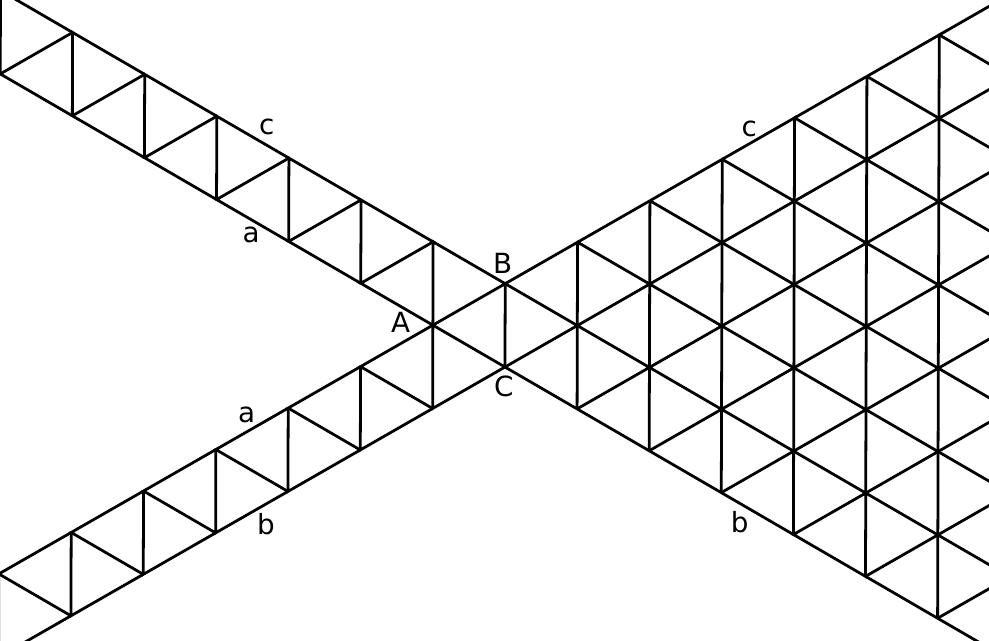}
\caption{\small A graph in $\mathcal{PC}_{\geq0}$ with total curvature $\frac{10}{12}$.}
\label{10d12}
\end{center}
\end{figure}

In this section, we construct some examples to show the sharpness of Theorem~\ref{mainthm2}.
In particular, we construct semiplanar graphs with nonnegative curvature whose total curvatures attain all the values in $$\left\{\frac{i}{12}:\ 0\leq i\leq 12, i\in \Z\right\}.$$
It is easy to construct examples for $i=0,12$. For $i=1$ and $11,$
we refer to Figure~\ref{fig2} and Figure~\ref{11d12-1} in the introduction. Examples with total curvature $\frac{i}{12}$ for $2\leq i\leq 10$ are shown in Figures~\ref{2d12}-\ref{10d12}. Note that half lines with same labeling are identified to each other.

\section{Semiplanar graphs with boundary}\label{s:boundary}

In this section, we study the total curvature for semiplanar graphs with boundary. 
Let $(V,E)$ be a graph topologically embedded into a surface $S$ with boundary. We call $G=(V,E,F)$ the semiplanar graph with boundary, where $F$ is the set of faces induced by the embedding. In this paper, we only consider semiplanar graphs satisfying tessellation properties, i.e. (i), (iii) in the introduction, and the condition (ii) replaced by the following
\begin{enumerate}[(ii')] \item The boundary of $S,$ denoted by $\partial S$, consists of edges of the graph. Each edge in the grap, by removing two end-vertices, is contained either in $\partial S$ or in the interior of $S$, denoted by $\mathrm{int}(S).$ It is incident to one face in the first case and to two different faces in the second.
\end{enumerate} Moreover, we always assume that $2\leq \deg(x)<\infty$ for any $x\in V\cap \partial S,$
$3\leq \deg(x)<\infty$ for any $x\in V\cap \mathrm{int}(S),$ and $3\leq \deg(\sigma)<\infty,$ for $\sigma\in F.$ As in the case of semiplanar graphs without boundary, one can define the  polygonal surface $S(G)$ for a semiplanar graph $G$ with boundary.


Let $G$ be a semiplanar graph with boundary, we define the combinatorial curvature as follows:
\begin{itemize}\item For any vertex $x\in \mathrm{int}(S),$ $\Phi(x)$ is defined as in \eqref{def:comb}.
\item  For any vertex $x$ on $\partial S$, 
\begin{equation}\label{def:curvbound}\Phi(x)=1-\frac{\deg(x)}{2}+\sum_{\sigma\in F:x\in \overline{\sigma}}\frac{1}{\mathrm{deg}(\sigma)}.\end{equation} 
\end{itemize} Direct calculation shows that for the vertex $x$ on the boundary of $S,$
\begin{equation}\label{eq:gaubonboundary}\Phi(x)=\frac{\pi-\theta_x}{2\pi},\end{equation}
where $\theta_x$ is the inner angle at $x$ w.r.t. $S$, as in \eqref{eq:GaussBonnetboundary}.

Adopting the same argument in \cite[Corollary 3.4]{HJL15}, one can show that a semiplanar graph with boundary $G$ has nonnegative combinatorial curvature if and only if the polygonal surface $S(G)$ has nonnegative sectional curvature in the sense of Alexandrov, see \cite{BuragoGromovPerelman92,BuragoBuragoIvanov01}. Let $G$ be a semiplanar graph with boundary and with nonnegative combinatorial curvature. Consider the doubling constructions of $S(G)$ and $G,$ see e.g. \cite{Perelman91}. Let $\widetilde{S(G)}$ be the
double of ${S}(G)$, that is, $\widetilde{S(G)}$ consists of two copies of ${S}(G)$ glued along the boundaries, which induces the doubling graph of $G,$ denoted by $\widetilde{G}.$ One can prove that $\widetilde{G}$ has nonnegative curvature, which can be derived from either the definition of the curvature on the boundary \eqref{def:curvbound} or Perelman's doubling theorem for Alexandrov spaces \cite{Perelman91}, and $\Phi(\widetilde{G})=2\Phi(G).$
This yields that $\widetilde{G}$ is a planar graph with nonnegative combinatorial curvature if the total curvature of $G$ is positive, which yields that $S(G)$ is homeomorphic to a half-plane with boundary. For our purposes, we always assume that $S$ is homeomorphic to a half-plane with boundary. Note that the double graph $\widetilde{G}$ may have vertices of degree two if $G$ has some vertices of degree two on $\partial S.$ However, the Cohn-Vossen type theorem, \cite[Theorem~1.3]{DeVosMohar07}, still applies to that case and implies that ${\Phi}(\widetilde{G})\leq 1.$ This yields 
\begin{equation}\label{eq:equ1}\Phi(G)\leq \frac12.\end{equation}

The patterns of a vertex on $\partial S$ with nonnegative curvature are given in the following proposition.
\begin{prop}\label{prop:curvatureboudary} Let $G$ be a semiplanar graph with boundary. If a vertex $x\in\partial S$ such that $\deg(x)\geq3$ and $\Phi(x)\geq0$, then the possible patterns of $x$ are
\begin{equation}\label{dx=2}
(3,3), (3,4), (3,5), (3,6), (4,4), (3,3,3),
\end{equation}
where the corresponding curvatures are
\begin{equation}\label{c_bd}
\frac{1}{6}, \frac{1}{12}, \frac{1}{30}, 0,0,0.
\end{equation}
\end{prop}

In the next lemma, we deal with the case that there is a vertex on $\partial S$ of degree two.
\begin{lemma}\label{lem:deg2}
Let $G$ be an infinite semiplanar graph with boundary and with nonnegative curvature, and $x$ be a vertex on $\partial S$ with $\deg(x)=2.$ Then the face, which $x$ is incident to, is of degree at most $6$ and $$\Phi(x)\geq\frac{1}{6}.$$ \end{lemma}
\begin{proof}
Let $x$ be incident to a $k$-gon $\sigma$. It suffices to prove that $k\leq 6.$ 

We claim that there is a vertex $y\in \partial \sigma\cap \partial S$ of degree at least three. Suppose it is not ture, i.e. $\deg(y)=2$ for any vertex $y\in \partial \sigma\cap \partial S,$ then we can show that \begin{equation}\label{eq:eq111}\partial \sigma\subset \partial S.\end{equation} In fact, we have the following observation: Given a vertex $y\in \partial \sigma\cap \partial S$ with $\deg(y)=2,$ any neighbor of $y,$ say $z,$ is in $\partial \sigma\cap \partial S$ and the edge connecting $y$ and $z$ is also contained in $\partial \sigma\cap \partial S.$ Then \eqref{eq:eq111} follows from the connectedness of $\partial \sigma$ and $\deg(y)=2$ for all vertices $y\in \partial \sigma\cap \partial S.$ This yields that $\partial S=\partial\sigma$ and $S=\sigma.$ This contradicts to the assumption that $G$ is infinite. This proves the claim.

Applying Proposition~\ref{prop:curvatureboudary} at a vertex on $\partial \sigma$ of degree at least three, we have $k\leq 6$ and hence prove the lemma.
\end{proof}

This lemma yields the following corollary.
\begin{corollary}\label{coro:maxdegree} Let $G=(V,E,F)$ be an infinite semiplanar graph with boundary and with nonnegative curvature. Then $$\deg(\sigma)\leq 42,\quad\forall\ \sigma\in F.$$
\end{corollary}
\begin{proof} Suppose not, i.e. there exists a face $\sigma$ such that $\deg(\sigma)\geq 43,$ then by Proposition~\ref{prop:curvatureboudary} and Lemma~\ref{lem:deg2}, the closure of $\sigma$ is contained in $\mathrm{int}(S).$ This yields $\Phi(G)\geq 1$ by Lemma~\ref{lemma}. It contradicts to \eqref{eq:equ1} and we prove the corollary.
\end{proof}

The next lemma shows that there are only finitely many vertices with non-vanishing curvature for a semiplanar graph with boundary and with nonnegative curvature, analogous to \cite[Theorem~1.4]{ChenChen08} and \cite[Theorem~3.5]{Chenbf09}. For a semiplanar graph $G$ with boundary, we denote by $T(G)$ the set of vertices with non-vanshing curvature.
\begin{lemma} Let $G$ be an infinite semiplanar graph with boundary and with nonnegative curvature. Then $T(G)$ is a finite set.
\end{lemma}
\begin{proof} By Corollary~\ref{coro:maxdegree}, the maximum facial degree is at most $42.$ For any vertex $x\in T(G),$ $\Phi(x)\geq \frac{1}{1722}$ if $x\in \mathrm{int}(S)$ by Table~\ref{tabl1}; $\Phi(x)\geq \frac{1}{30}$ if $x\in \partial S$ by Proposition~\ref{prop:curvatureboudary} and Lemma~\ref{lem:deg2}. Denote by $N$ the number of vertices in $T(G).$ Then
$$\frac{N}{1722}\leq \sum_{x\in T(G)}\Phi(x)=\Phi(G)\leq \frac{1}{2}.$$ This proves that $N\leq 861.$
\end{proof}

Now we are ready to prove the main results in the section. We determine all possible values of total curvatures of semiplanar graphs with boundary and with nonnegative curvature in the next theorem.
\begin{theorem}\label{half-plane-2}
Let $G$ be an infinite semiplanar graph with boundary and with nonnegative curvature. 
Then the total curvature of $G$ is an integral multiple of $\frac{1}{12}.$ 
\end{theorem}

\begin{proof} We use the similar arguments as in the proof of Theorem \ref{mainthm2}. Without loss of generality, we assume that $S(G)$ is homeomorphic to the half-plane with boundary. One can find a (sufficient large) compact set $K\subset S(G)$ satisfying the following:
\begin{itemize}
\item $K$ is the closure of a finite union of faces in $S(G)$ and is simply connected.
\item $T(G)\subset K.$
\item $\partial K$ is a Jordan curve consisting of finitely many edges.
\item $\partial K\cap \partial S$ contains at least one edge and is connected. 
\item For any vertex $x\in \overline{\partial K\setminus \partial S},$ $B_2(x)\cap T(G)=\emptyset.$
\end{itemize} 
Set $V_1(K):=V\cap \overline{\partial K\setminus \partial S}$ and $V_2(K):=(V\cap K)\setminus V_1(K).$ The Gauss-Bonnet formula on $K$ reads as
$$2\pi\sum_{x\in V_2(K)\setminus \partial S}\Phi(x)+\sum_{x\in V_2(K)\cap \partial S}(\pi-\theta_x)+\sum_{x\in V_1(K)}(\pi-\theta_x)=2\pi,$$ where $\theta_x$ is the inner angle of $x$ w.r.t. the set $K.$
By \eqref{eq:gaubonboundary}, we may rewrite the above equation as
$$2\pi\sum_{x\in V_2(K)}\Phi(x)+\sum_{x\in V_1(K)}(\pi-\theta_x)=2\pi.$$
By the properties of $K,$ $\sum_{x\in V_2(K)}\Phi(x)=\sum_{x\in T(G)}\Phi(x)=\Phi(G).$

To prove the theorem, it suffices to show that $\sum_{x\in V_1(K)}(\pi-\theta_x)$ is an integral multiple of $\frac{\pi}{6}$. For any $x\in V_1(K)\cap \partial S,$ since $\Phi(x)=0,$ Proposition~\ref{prop:curvatureboudary} implies that the pattern of $x$ is $(3,6), (4,4)$ or $(3,3,3).$ 
We claim that for any $x\in V_1(K)\setminus \partial S,$ the possible patterns of $x$ are as follows:
$$(3,12,12), (4,6,12), (6,6,6), (3,3,4,12), (3,3,6,6), (3,4,4,6),$$
$$(4,4,4,4), (3,3,3,3,6), (3,3,3,4,4), (3,3,3,3,3,3).$$ Let $x\in V_1(K)\setminus \partial S.$ Consider the double of $G,$ denoted by $\widetilde{G}.$ This induces the double of $K,$ $\widetilde{K},$ whose boundary is a Jordan curve. One can show that for any vertex $y\in \partial \widetilde{K},$ the vertices in $\widetilde{B}_2(y),$ the ball of radius two centered at $y$ in the double $\widetilde{G},$ have vanishing curvature and are of vertex degree at least three. By Lemma~\ref{lemma2.1}, we have the pattern list for $y$ therein (and hence for $x$). To prove the claim, it suffices to exclude the pattern $(4,8,8)$ for the vertex $x.$ We argue by contradiction. Suppose that $x=(4,8,8),$ then by Lemma~\ref{lem:488}, $y=(4,8,8)$ for any $y\in \partial \widetilde{K}.$  Since $\partial K\cap \partial S$ contains at least one edge, and hence, there is a vertex $y_0\in \partial \widetilde{K}\cap \partial S.$ However, the pattern of $y_0$ in $G$ is $(3,6), (4,4)$ or $(3,3,3)$ which yields a contradiction. This proves the claim.

Hence, as in Case 2 in the proof of Theorem~\ref{mainthm2}, by dividing dodecagons and hexagons into triangles and squares, we may assume that any face, which some vertex in $V_1(K)$ is incident to, is either a triangle or a square. Thus the inner angle of any vertex in $V_1(K)$ is of the form $m\frac{\pi}{2}+n\frac{\pi}{3}$ $(m,n\in\mathbb{Z}),$ which yields that $\sum_{x\in V_1(K)}(\pi-\theta_x)$ is an integral multiple of $\frac{\pi}{6}$. This proves the theorem.

\end{proof}

\begin{figure}[htbp]
\begin{center}
\includegraphics[width=0.5\linewidth]{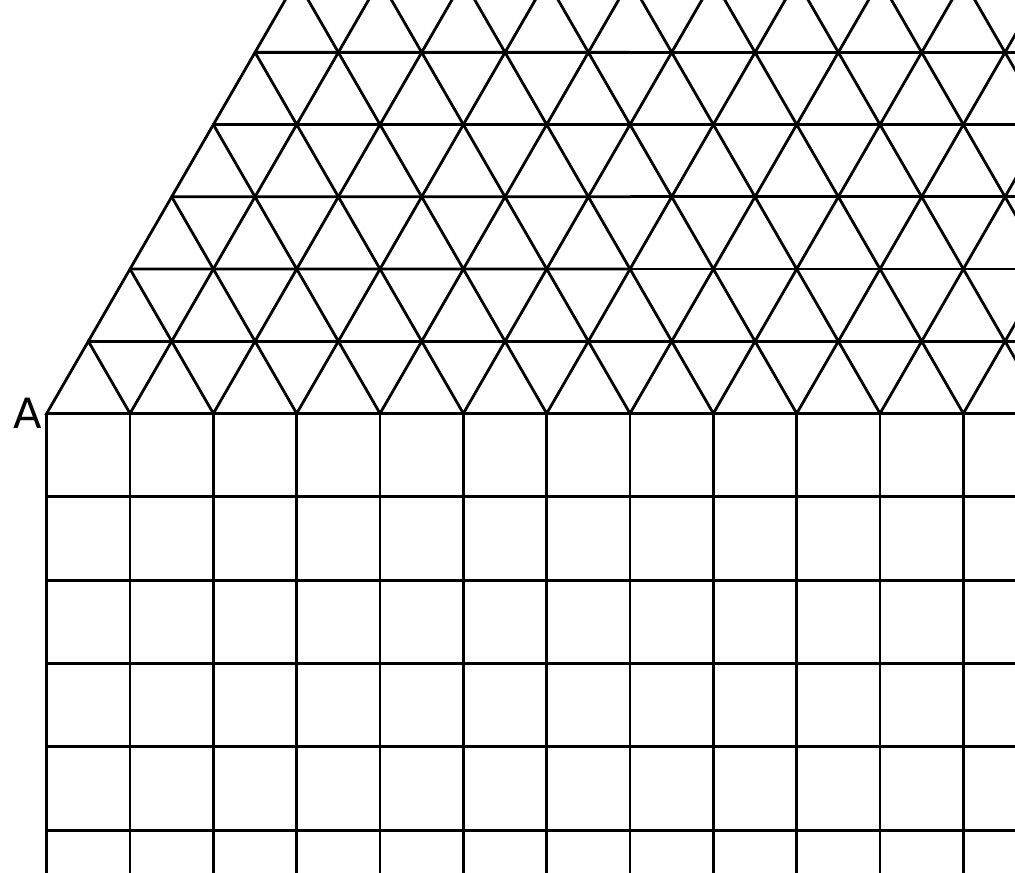}
\caption{\small A graph with boundary and with total curvature $\frac{1}{12}$.}
\label{fig-half1}
\end{center}
\end{figure}

\begin{figure}[htbp]
\begin{center}
\includegraphics[width=0.6\linewidth]{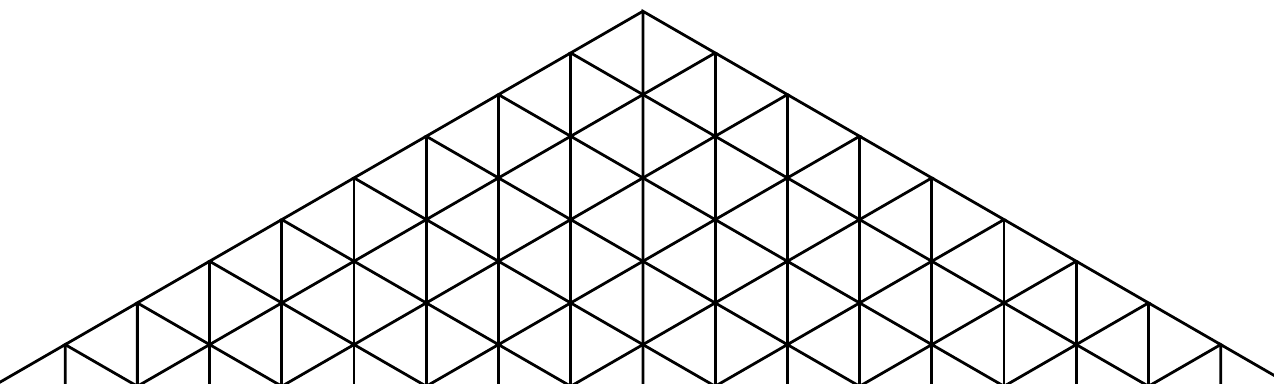}
\caption{\small A graph with boundary and with total curvature $\frac{2}{12}$.}
\label{fig-half2}
\end{center}
\end{figure}

\begin{figure}[htbp]
\begin{center}
\includegraphics[width=0.55\linewidth]{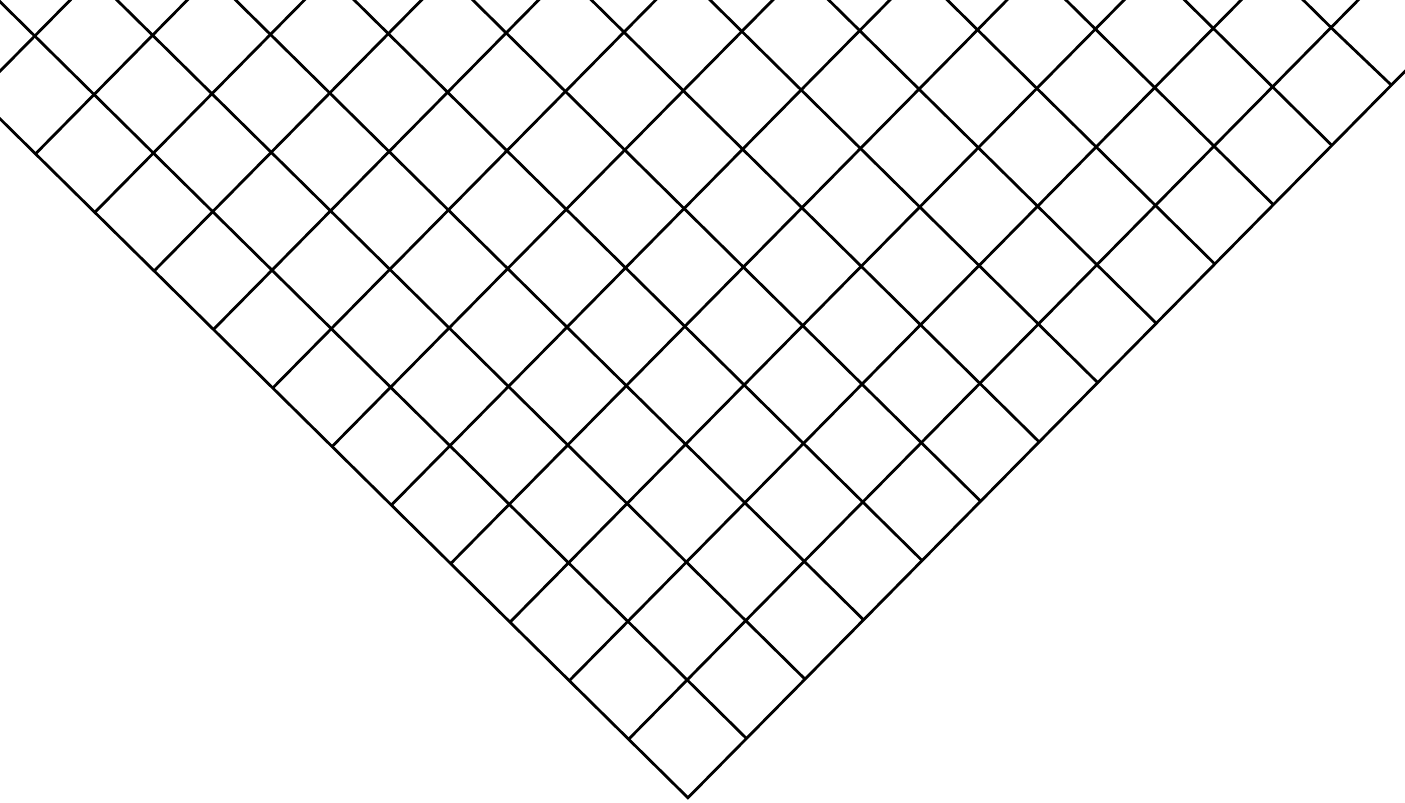}
\caption{\small A graph with boundary and with total curvature $\frac{3}{12}$.}
\label{fig-half3}
\end{center}
\end{figure}

\begin{figure}[htbp]
\begin{center}
\includegraphics[width=0.5\linewidth]{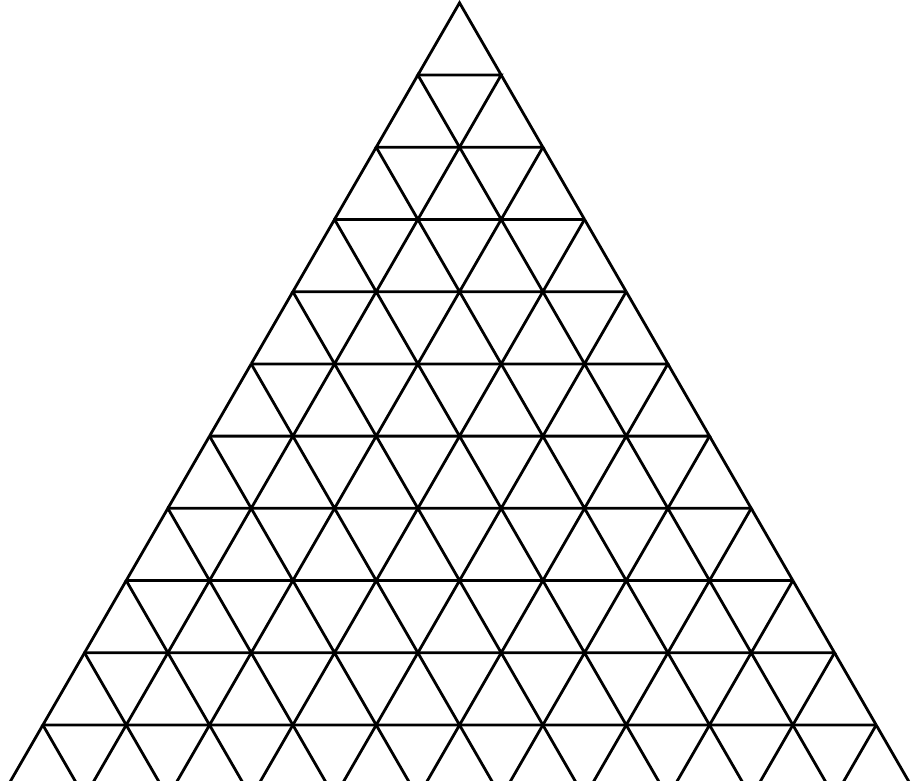}
\caption{\small A graph with boundary and with total curvature $\frac{4}{12}$.}
\label{fig-half4}
\end{center}
\end{figure}

\begin{figure}[htbp]
\begin{center}
\includegraphics[width=0.5\linewidth]{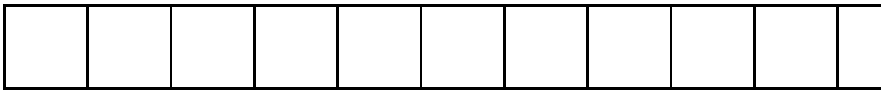}
\caption{\small A graph with boundary and with total curvature $\frac{6}{12}$.}
\label{fig-half6}
\end{center}
\end{figure}

By this theorem, the total curvature of a semiplanar graph with boundary and with nonnegative curvature takes the value in $$\left\{\frac{i}{12}:\ 0\leq i\leq 6, i\in \Z\right\}.$$ It is easy to construct an example for $i=0.$ We construct examples for $i=1,2,3,4,6,$ see Figures~\ref{fig-half1}-\ref{fig-half6}. We ask whether one can construct an example with total curvature $\frac{5}{12}$ and propose the following conjecture.
\begin{conjecture} There is no semiplanar graph with boundary and with nonnegative curvature whose total curvature is $\frac{5}{12}.$
\end{conjecture}


At the end of the section, we classify the graph/tessellation structures of semiplanar graphs with boundary and with nonnegative curvature whose total curvatures attain the first gap of the total curvature, i.e. $\frac{1}{12}$. 

\begin{figure}[htbp]
\begin{center}
\includegraphics[width=0.65\linewidth]{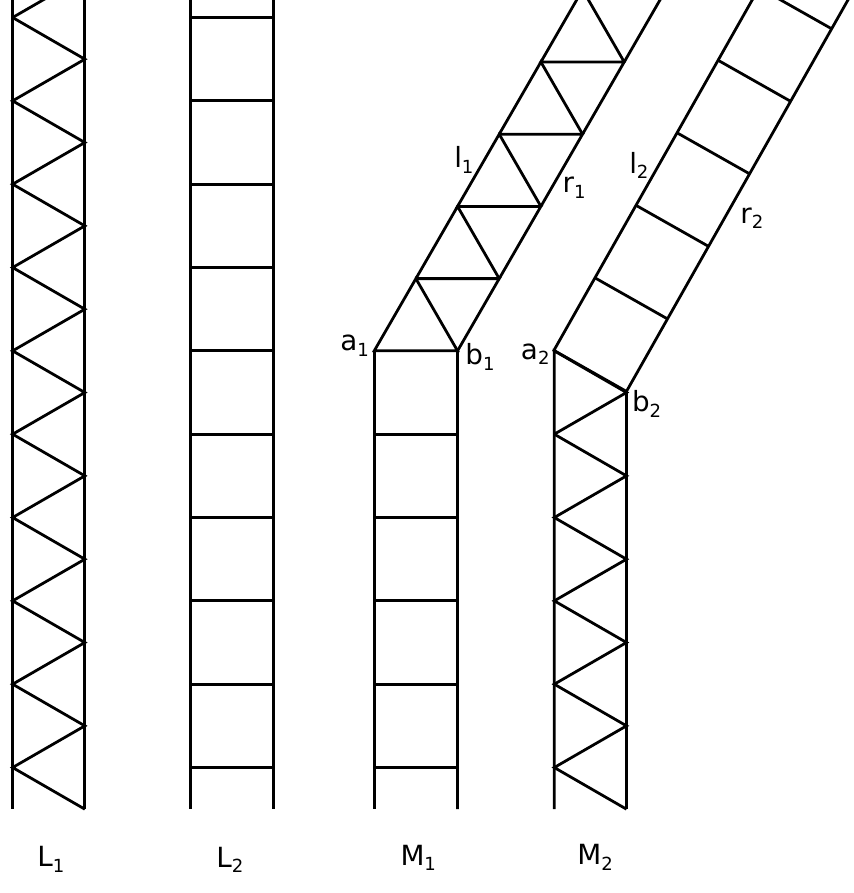}
\caption{\small The strips: $L_1$, $L_2$, $M_1$ and $M_2$.}
\label{L-M}
\end{center}
\end{figure}

\begin{figure}[htbp]
\begin{center}
\includegraphics[width=0.7\linewidth]{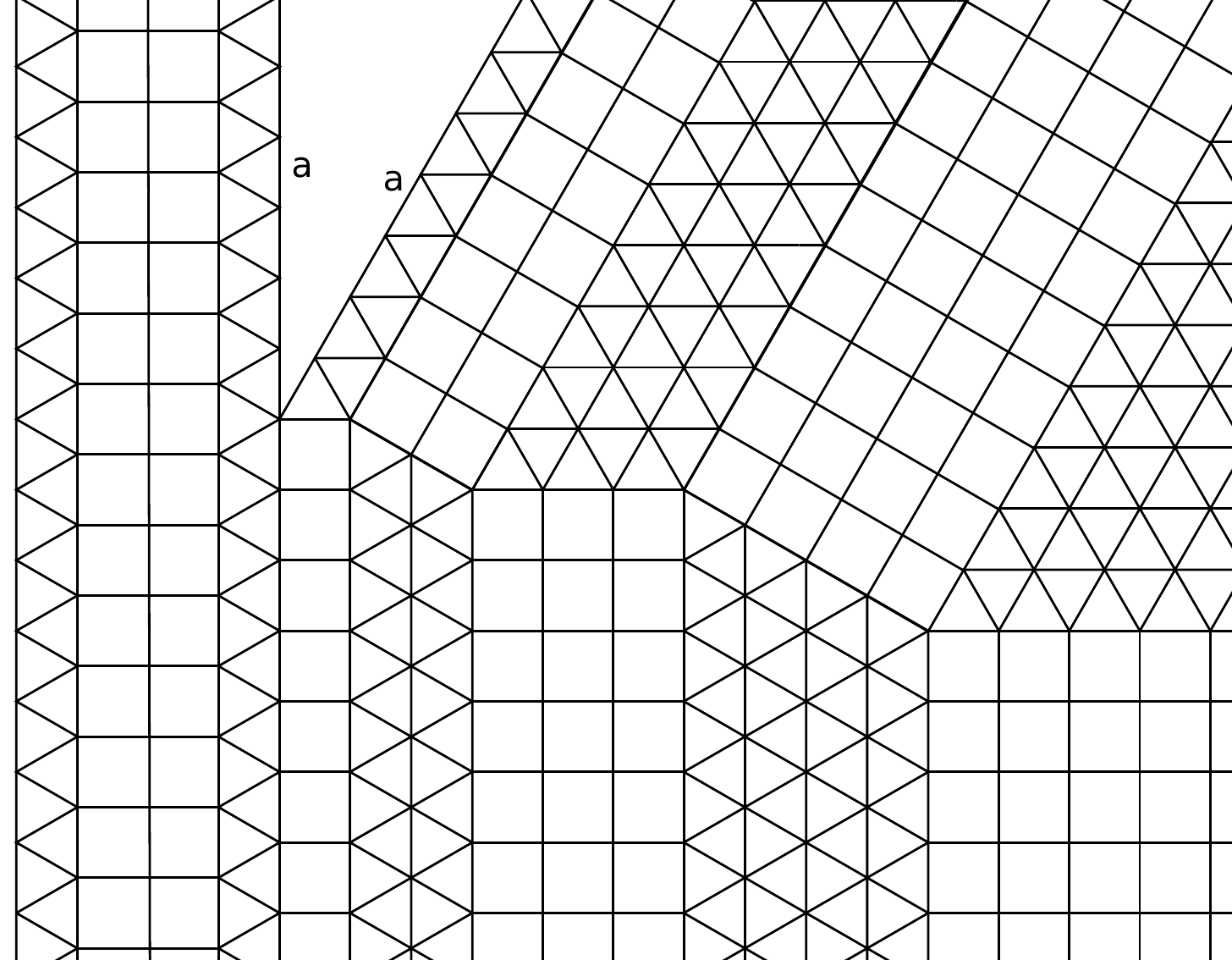}
\caption{\small An example obtained by gluing a sequence of $L_1,$ $L_2,$ $M_1$ and $M_2.$}
\label{gluingsequence}
\end{center}
\end{figure}

We denote by $L_1$ ($L_2$ resp.) the strip consists of triangles (squares resp.), and by $M_1$ ($M_2$ resp.) the strip with triangles in the upper half and squares in the lower half (with squares in the upper half and triangles in the lower half resp.), see Figure \ref{L-M}. 
We denote by $l_1$ and $r_1$ ($l_2$ and $r_2$ resp.) the left part and the right part of the boundary of $M_1$ ($M_2$ resp.), and by $a_1$ and $b_1$ ($a_2$ and $b_2$ resp.) the vertices on $l_1$ and on $r_1$ ($l_2$ and $r_2$ resp.) incident to one or two triangles and a square. 
We define the following gluing rules for these four strips:
\begin{itemize} \item For any strip $L_\alpha,$ $\alpha=1$ or $2,$ it can be glued to $L_\beta$ or $M_\gamma$ for $\beta=1,2$ and $\gamma=1,2,$ from the left to the right, denoted by $L_{\alpha}L_\beta$ or $L_{\alpha}M_\gamma.$ 
\item For any strip $M_\alpha,$  $\alpha=1$ or $2,$ it can be glued to $M_\gamma$ for $\gamma=1,2,$ from the left to the right, such that $r_\alpha$ is identified with $l_\gamma$ and $b_\alpha$ is identified with $a_\gamma,$ denoted by $M_{\alpha}M_\gamma.$
\end{itemize} An example obtained by gluing a sequence $$L_1L_2L_2L_1M_1M_2M_2M_1M_1M_1M_2M_2M_2M_2M_1\cdots,$$ is shown in Figure~\ref{gluingsequence}.

For any semiplanar graph $G$ and the polygonal surface $S(G),$ for each hexagon in $G$ we add a new vertex at the barycenter and join it to the vertices of the hexagon to obtain a semiplanar graph, denoted by $\widehat{G}.$ In this way, $\widehat{G}$ is obtained from $G$ by replacing each hexagon in $G$ with six triangles.
\begin{theorem}\label{thm:classification}
Let $G$ be an infinite semiplanar graph with boundary and with nonnegative curvature. Then the total curvature of $G$ is $\frac{1}{12}$ if and only if $\widehat{G}$ is isomorphic to one of the following:
\begin{enumerate}
\item $P_1P_2\cdots P_n\cdots,$ 
\item $Q_1\cdots Q_TP_1P_2\cdots P_n\cdots,$ for some $T\in\N.$
\end{enumerate} where $Q_t=L_1$ or $L_2$ for any $1\leq t\leq T,$ $P_n=M_1$ or $M_2$ for any $n\in \N,$ and the successive strips are glued together by the way defined before.

\end{theorem}
\begin{proof} The ``if" part of the theorem is obvious. We only need to show the ``only if" part.
It suffices to consider semiplanar graphs with no hexagonal faces. By Lemma~\ref{lem:deg2}, $\Phi(G)=\frac{1}{12}$ implies that all vertices of $G$ have vertex degree at least three. We divide it into two cases.

\begin{description}
\item[Case 1] There exists a vertex $A$ on $\partial S$ with non-vanishing curvature. Then the possible patterns of $A$ are $(3,3), (3,4)$ and $(3,5)$. We claim that $A=(3,4).$ If $A=(3,3)$, then $\Phi(A)=\frac{1}{6}>\frac{1}{12},$ which contradicts to $\Phi(G)=\frac{1}{12}$. 

\begin{figure}[tb]
\begin{minipage}[b]{0.49\textwidth}
\centering
\includegraphics[width=2cm,height=3cm]{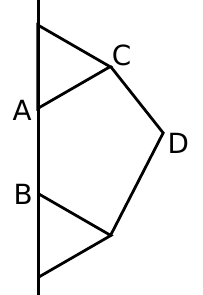}
\end{minipage}
\caption{\small Proof of Theorem~\ref{thm:classification}}
 \label{ex35}
\end{figure}

If $A=(3,5)$, see Figure \ref{ex35}, then $B=(3,5)$ and $\Phi(A)=\Phi(B)=\frac{1}{30}.$ The possible patterns of $C$ are $(3,3,5)$, $(3,4,5)$, $(3,5,5)$, $(3,5,6)$, $(3,3,3,5)$, $(3,3,4,5)$, $(3,3,5,5)$, $(3,3,5,6)$, $(3,3,5,7)$, $(3,4,4,5)$ and $(3,3,3,3,5)$. However, except the case of $C=(3,3,5,7)$, we always have $\Phi(C)>\frac{1}{60}$. Hence $\Phi(G)\geq\Phi(A)+\Phi(B)+\Phi(C)>\frac{1}{12}$. We only need to consider $C=(3,3,5,7).$ In this case, $\Phi(C)=\frac{1}{105}.$ Moreover, one can show that $\Phi(D)\geq\frac{1}{105}$, which yields $\Phi(G)\geq\Phi(A)+\Phi(B)+\Phi(C)+\Phi(D)>\frac{1}{12}$. This also yields a contradiction and hence prove the claim. 

We consider the case of $A=(3,4).$ In this case, $\Phi(A)=\frac{1}{12}$ and hence any other vertex has vanishing curvature. To visualize the graph, we draw the picture, from the left to the right, on the plane. Let two faces, a triangle and a square which $A$ is incident to, locate on the right half plane passing through $A,$ such that two edges on $\partial S,$ which $A$ is incident to, are in the direction of $-y$ axis (i.e. $\theta=-\frac{\pi}{2}$) and $\theta=\frac{\pi}{3}$ respectively, where $\theta$ is the angular coordinate in polar coordinates. We have two 
subcases: 
\begin{enumerate}[(a)]\item The triangle is above the square.
\item The square is above the triangle.
\end{enumerate} Without loss of generality, we consider the subcase $(a).$ Denote by $X_1$ ($Y_1$ resp.) a neighbor of $A$ on $\partial S$ which is incident to the triangle (the square resp.). Since any vertex except $A$ has vanishing curvature, one can show that $X_1=(3,3,3)$ ($Y_1=(4,4)$ resp.). Set $X_0=Y_0=A.$ For any $i\in\N,$ we inductively define $X_{i+1}$ ($Y_{i+1}$ resp.) as the neighbor of $X_i$ ($Y_{i}$ resp.) which is different from $X_{i-1}$ ($Y_{i-1}$ resp.) By the induction, one can show that $X_i=(3,3,3)$ and $Y_i=(4,4)$ for $i\in \N.$ This yields the strip $M_1,$ as in Figure \ref{L-M}. Similarly, in the subcase $(b)$ we have the strip $M_2.$ We cut off the strip $M_1$ (or $M_2$) from the semiplanar graph $G$ to obtain a new graph $G_1=(V_1,E_1,F_1),$ precisely,
\begin{eqnarray*}&&V_1=V\setminus \partial S, \quad E_1=\{\{x,y\}\in E: x,y\in \mathrm{int}(S)\},\\ &&F_1=\{\sigma\in F: \sigma\subset \mathrm{int}(S)\}.\end{eqnarray*}
One can check that $G_1$ has nonnegative curvature, the total curvature of $G_1$ is $\frac{1}{12},$ and the neighbor of $A$ in $G,$ which is not $X_1$ and $Y_1,$ is of curvature $\frac{1}{12}.$
Hence, this reduces the case of $G_1$ to Case~1 with $G=G_1$, and we can cut off a strip, $M_1$ or $M_2,$ from $G_1$ to get a new graph $G_2.$ Continuing this process by cutting off strips along the boundary, we get $(1)$ in the theorem.

\item[Case 2] All vertices on $\partial S$ have vanishing curvature.  In this case, the patterns of vertices along the boundary are $(3,3,3)$ or $(4,4)$. If there exists a vertex on the boundary of pattern $(3,3,3)$ ((4,4) resp.), then the patterns of other vertices along the boundary are also $(3,3,3)$ ((4,4) resp.). This yields the strip $L_1$ ($L_2$ resp.) along the boundary. As in Case 1, we cut off the strip $L_1$ (or $L_2$) from the graph to obtain a new graph, which has nonnegative curvature and the total curvature $\frac{1}{12}.$ For the new graph, if all vertices on the boundary have vanishing curvature, then we are in the situation of Case 2. So that we may cut off the strip $L_1$ (or $L_2$) in the new graph to get another graph. We continue this process unless we find some vertex on the boundary with non-vanishing curvature. This process stops in finite steps since the total curvature is positive and there is a vertex of non-vanishing curvature. In the last step, we are in the situation of Case 1 and get the classification accordingly. By the re-construction, we have $(2)$ in the theorem.
\end{description}
\end{proof}


\bigskip
{\bf Acknowledgements.} We thank Tam\'as R\'eti for many helpful discussions on the problems.

B. H. is supported by NSFC, grant no. 11401106. Y. S. is supported by NSF of Fujian Province through Grants 2017J01556, 2016J01013, JAT160076.

\bibliography{Reti-ref}
\bibliographystyle{alpha}

\end{document}